\documentstyle[12pt]{article}
\makeatletter \oddsidemargin  -.1in \evensidemargin -.1in
\newcommand{\singlespacing}{\let\CS=\@currsize\renewcommand{\baselinestretch}{1}\tiny\CS}
\newcommand{\oneandahalfspacing}{\let\CS=\@currsize\renewcommand{\baselinestretch}{1.25}\tiny\CS}
\newcommand{\doublespacing}{\let\CS=\@currsize\renewcommand{\baselinestretch}{1.35}\tiny\CS}

\newtheorem{rule-def}[theorem]{Rule}

\RequirePackage[dvips]{graphicx} \textheight 21.5cm

\begin{document}


\title{\bf Pulsatile flow and heat transfer of a magneto-micropolar fluid through a stenosed
artery under the influence of body acceleration}
\author{\small G. C. Shit\thanks{Corresponding author. Email address: gcs@math.jdvu.ac.in
(G. C. Shit)}  ~~and ~~M. Roy\\
 \it Department of Mathematics,\\ Jadavpur University, Kolkata
- 700032, India
\\}
\date{}
\maketitle \noindent \doublespacing
\begin{abstract} With an aim to investigate the effect of externally
imposed body acceleration and magnetic field on pulsatile flow of
blood through an arterial segment having stenosis is under
consideration in this paper. The flow of blood is presented by a
unsteady micropolar fluid and the heat transfer characteristics
have been taken into account. The non-linear equations that
governing the flow are solved numerically using finite difference
technique by employing a suitable coordinate transformation. The
numerical results have been observed for axial and microrotation
component of velocity, fluid acceleration, wall shear stress(WSS),
flow resistance, temperature and the volumetric flow rate. It thus
turns out that the rate of heat transfer increases with the
increase of Hartmann number $H$, while the wall shear stress has a
reducing effect on the Hartmann number $H$ and an enhancing effect
on microrotation parameter $K$ as well as the
constriction height $\delta$.\\

\end{abstract}

\section {Introduction:}
It is well known that fluid dynamical behaviour of blood through
an arterial segment having stenosis play a vital role in
cardiovascular disease. The narrowing in the artery, commonly
referred to as stenosis, is a dangerous disease and is caused due
to the abnormal growth in the lumen of the arterial wall. Stenosis
may be formed at one or more locations of the cardiovascular
system. As a result of such undesirable formation at the
endothelium of the vessel wall, reduction of regular blood flow is
likely to take place near the stenosis. If this disease takes a
severe form, it may lead to stroke, heart attack and various
cardiovascular disease. It has been observed that whole blood,
being a suspension of blood cells in plasma, behave as a
non-Newtonian fluid. Several attempts [1-8] have been made to
understand the flow characteristic of blood through arteries under
various assumptions. But most of these studies have failed to
estimate on an important role in blood rheology is the motion of
erythrocytes, white blood cells
and platelets in plasma.\\
Eringen \cite{R10,R11} introduced the theory of micropolar fluid,
which exhibits microscopic effects arising from the local
structure and microrotation of fluid microelements are considered
to be rigid. A subclass of microfluids such as liquid crystals,
suspensions and animal blood are treated as micropolar fluid,
which takes into account the rotation of fluid particles by means
of an independent kinematic vector called the microrotation
vector. With this end in view, several investigators [12-14]
considered micropolar fluid as a blood flowing through an arterial
stenosis. They suggested that blood exhibits remarkable
non-Newtonian properties, since blood is a suspension of neutrally
bouyant deformable particles in viscous fluid. Devanathan and
Parvathamma \cite{R16} developed a mathematical model for steady
flow of micropolar fluid through an arterial stenosis and their
results indicated that the axial velocity is highly affected owing
to the presence of microelements in the fluid. However, all the
above studies did not considered the effect of body acceleration.
Although, under normal physiological conditions, flow properties
of blood depends upon pumping action of the heart, which in turn
produces a pulsatile pressure gradient through the cardiovascular
system. Although, there may be some another important situation
under which human beings experiences whole body acceleration. It
may occur for a short or a long period of time. Such situations
arise day to day life of human being e.g., driving a vehicle,
flying in an aircraft, due to which there may occur, the serious
health problems like headache, loss of vision, increase of pulse
rate and hemarage in face, neck and eye socket. Keeping this in
mind, several investigators [16-24] carried out their studies
pertaining to blood flow through stenosed arteries with a periodic
body acceleration. Chakravarty and Sanniggrahi \cite{R27} studied
theoretically the flow characteristics of blood through an artery
in the presence of multi-stenosis when it is subjected to whole
body acceleration. Wherein they considered a flexible cylindrical
tube containing a homogeneous Newtonian fluid with variable
viscosity characterising blood. But all the above studies neither
consider the effect of body acceleration on pulsating blood flow
nor includes the effect of magnetic field when blood is considered to
be micropolar fluid model.\\
Since blood consists of suspension of red cells containing
hemoglobin, which contains iron oxide, it is quite apparent that
blood is electrically conducting and exhibits magnetohydrodynamic
flow characteristics. If a magnetic field is applied to a moving
and electrically conducting fluid, it will induce electric as well
as magnetic fields. The interaction of these fields produce a body
force  known as Lorentz force, which has a tendency to slow down
the fluid motion. Such analysis may be useful for pumping of blood
and magnetic resonance imaging (MRI). The magnetic field can also
be used for controlling blood flow during surgery. Many authors
[26-29] have investigated the flow of blood through arteries in
the presence of magnetic field under different situations. The
steady MHD flow of a viscous fluid in a slowly varying channel in
the presence of a uniform magnetic field was explored by Misra et
al \cite{R33}. Recently Misra and Shit \cite{R34} investigated the
flow and heat transfer of a viscoelastic electrically conducting
fluid under the action of a magnetic field. They made an
observation that the temperature of blood increases as the
strength of magnetic field increases is likely to be useful in the
development of new heating methods. \\
With the above discussion in mind, a good attempt is made in the
present theoretical investigation to examine the pulsatile blood
flow characteristics through a stenosed artery experiencing
periodic body acceleration by treating micropolar fluid model for
blood in the presence of magnetic field. The motivation of this
study is to analyze the effect of microrotation component and
axial velocity when an external magnetic field and body
accelerations are applied. In addition, in view of the information
stated above, we have also performed a heat transfer phenomena of
the problem in question.

\section{Mathematical Formulation}
Let us consider a fully developed unsteady, laminar,
incompressible and axially symmetric two dimensional flow through
a stenosed artery. An external magnetic field is applied to a
pulsatile flow of blood, which is characterized by a micropolar
fluid model under periodic body acceleration. Let $(r^{*},
\theta^{*}, z^{*})$ be the co-ordinates of the material point in a
cylindrical polar co-ordinate system, of which $z^{*}$ is taken as
central axis of the artery. We denote $u^{*}$ and $v^{*}$ the
velocity components of blood along the axial and radial directions
respectively and $w^{*}$ the component of microrotation.  For
axisymmetry, we assumed that the flow variation is independent of
$\theta^{*}$. Owing to the above considerations, the governing
equations of the present problem may be written in the form,
\begin{eqnarray}
         \frac{\partial u^{*}}{\partial
         z^{*}}+\frac{v^{*}}{r^{*}}+\frac{\partial v^{*}}{\partial
         r^{*}}=0
\end{eqnarray}
\begin{eqnarray}
     \rho(\frac{\partial u^{*}}{\partial t^{*}}+u^{*}\frac{\partial u^{*}}{\partial z^{*}}+v^{*}\frac{\partial u^{*}}{\partial
         r^{*}})=-\frac{\partial p^{*}}{\partial z^{*}}+(\mu+k)(\frac{\partial^{2} u^{*}}{\partial {r^{*}}^{2}}
         +\frac{1}{r^{*}}\frac{\partial u^{*}}{\partial r^{*}}+\frac{\partial^{2} u^{*}}{\partial {z^{*}}^{2}})+k(\frac{\partial w^{*}}{\partial
         r^{*}}+\frac{w^{*}}{r^{*}})\nonumber\\-\sigma B^{2}_{0}u^{*}+\rho
         G^{*}(t^{*}),
 \end{eqnarray}
\begin{eqnarray}
  \rho(\frac{\partial v^{*}}{\partial t^{*}}+u^{*}\frac{\partial v^{*}}{\partial z^{*}}+v^{*}\frac{\partial v^{*}}{\partial
         r^{*}})=-\frac{\partial p^{*}}{\partial r^{*}}+(\mu+k)(\frac{\partial^{2} v^{*}}{\partial {r^{*}}^{2}}
         +\frac{1}{r^{*}}\frac{\partial v^{*}}{\partial r^{*}}-\frac{v^{*}}{{r^{*}}^{2}}+\frac{\partial^{2} v^{*}}{\partial {z^{*}}^{2}})-k\frac{\partial w^{*}}{\partial
         z^{*}},
\end{eqnarray}
\begin{eqnarray}
   \rho j(\frac{\partial w^{*}}{\partial t^{*}}+u^{*}\frac{\partial w^{*}}{\partial z^{*}}+v^{*}\frac{\partial w^{*}}{\partial
         r^{*}})=-k(2w^{*}+\frac{\partial u^{*}}{\partial r^{*}}-\frac{\partial v^{*}}{\partial
         z^{*}})+\gamma(\frac{\partial^{2} w^{*}}{\partial {r^{*}}^{2}}+\frac{1}{r^{*}}\frac{\partial w^{*}}{\partial r^{*}}-\frac{w^{*}}{{r^{*}}^{2}}+\frac{\partial^{2} w^{*}}{\partial
         {z^{*}}^{2}}),
\end{eqnarray}
\begin{eqnarray}
    \frac{\partial T^{*}}{\partial t^{*}}+u^{*}\frac{\partial T^{*}}{\partial z^{*}}+v^{*}\frac{\partial T^{*}}{\partial
         r^{*}}=\frac{\kappa_0}{\rho C_{p}}(\frac{\partial^{2} T^{*}}{\partial {r^{*}}^{2}}+\frac{1}{r^{*}}\frac{\partial T^{*}}{\partial r^{*}}+\frac{\partial^{2} T^{*}}{\partial
         {z^{*}}^{2}})+\frac{\sigma B^{2}_{0}}{\rho C_{p}}{u^{*}}^{2}
\end{eqnarray}
where $\rho$ represents the blood density, $p$ the blood pressure,
$\mu$ the viscosity of blood, $k$ the rotational viscosity,
$\sigma$ the electrical conductivity, $B_{0}$ the applied magnetic
field intensity, $\gamma$ the spin gradient viscosity is assumed
to be $\gamma=j(\mu+k/2)$, $j$ the micro-gyration constant,
$T^{*}$ is the temperature, $C_{p}$ the specific heat and
$\kappa_0$ denotes thermal conductivity. In the above, we have
neglected the effect of induced magnetic field due to low magnetic
Reynolds number ($R_{m}<<1$). For ~$t^{*}>0$, the flow is assumed
to have a periodic body acceleration $G^{*}(t^{*})$ appeared in
equation (2) has the expression of the form
\begin{eqnarray}
               G^{*}(t^{*})=a^{*}(cos\omega_{b}t^{*}+\phi_{g})
\end{eqnarray}
where $a^{*}$ denote the amplitude of body acceleration,
$\omega_{b}$ the
frequency of body acceleration and $\phi_{g}$ the phase difference.\\
Let us introduce the following non-dimensional variables to put
the above equations in dimensionless form:
 \[z=\frac{z^{*}}{R_{0}},~~r=\frac{r^{*}}{R_{0}},~~u=\frac{u^{*}}{\omega R_{0}},~~v=\frac{v^{*}}{\omega
 R_{0}},~~w=\frac{w^{*}}{\omega},~~p=\frac{p^{*}}{\mu
 \omega},~~t=\omega
 t^{*},~~J=\frac{j}{R^{2}_{0}},~~\theta=\frac{T^{*}-T_{\infty}}{T_{w}-T_{\infty}}.\]
Using these dimensionless variables the equations (1)-(5) becomes
\begin{eqnarray}
          \frac{\partial u}{\partial z}+\frac{v}{r}+\frac{\partial v}{\partial r}=0
\end{eqnarray}
\begin{eqnarray}
   \alpha^{2}(\frac{\partial u}{\partial t}+u\frac{\partial u}{\partial z}+v\frac{\partial u}{\partial
         r})=-\frac{\partial p}{\partial z}+(1+K)(\frac{\partial^{2} u}{\partial {r}^{2}}
         +\frac{1}{r}\frac{\partial u}{\partial r}+\frac{\partial^{2} u}{\partial z^{2}})+K(\frac{\partial w}{\partial
         r}+\frac{w}{r})-H^{2}u+G(t)
\end{eqnarray}
\begin{eqnarray}
    \alpha^{2}(\frac{\partial v}{\partial t}+u\frac{\partial v}{\partial z}+v\frac{\partial v}{\partial
         r})=-\frac{\partial p}{\partial r}+(1+K)(\frac{\partial^{2} v}{\partial r^{2}}
         +\frac{1}{r}\frac{\partial v}{\partial r}-\frac{v}{r^{2}}+\frac{\partial^{2} v}{\partial z^{2}})-k\frac{\partial w}{\partial z}
\end{eqnarray}
\begin{eqnarray}
   \alpha^{2} J(\frac{\partial w}{\partial t}+u\frac{\partial w}{\partial z}+v\frac{\partial w}{\partial
         r})=-K(2w+\frac{\partial u}{\partial r}-\frac{\partial v}{\partial z})+m(\frac{\partial^{2} w}{\partial r^{2}}+\frac{1}{r}\frac{\partial w}{\partial r}-\frac{w}{r^{2}}+\frac{\partial^{2} w}{\partial z^{2}})
\end{eqnarray}

\begin{eqnarray}
  \frac{\partial \theta}{\partial t}+u\frac{\partial \theta}{\partial z}+v\frac{\partial \theta}{\partial r}
  =\frac{1}{\alpha^{2} P_{r}}(\frac{\partial^{2} \theta}{\partial r^{2}}+\frac{1}{r}\frac{\partial \theta}{\partial r}+\frac{\partial^{2} \theta}{\partial
         z^{2}})+\frac{E_{c}H^{2}}{\alpha^{2}}u^{2}
\end{eqnarray}\\
where the non-dimensional parameters appeared in equations
(8)-(11) are defined as the ratio of viscosity $K=\frac{k}{\mu}$,
the material constant $m=\frac{\gamma}{\mu R^{2}_{0}}$, Womersley
number $\alpha=R_{0}\sqrt{\frac{\rho \omega}{\mu}}$, Hartmann
number $H=R_{0}B_{0}\sqrt{\frac{\sigma}{\mu}}$, Prandtl
number $P_{r}=\frac{\mu C_{p}}{\kappa_{0}}$, and Eckert number $E_{c}=\frac{\omega^{2}R^{2}_{0}}{C_{p}(T_{w}-T_{\infty})}$.\\
Using the non-dimensional quantities $a_{0}=\frac{\rho
R_{0}a^{*}}{\mu \omega}$, $b=\frac{\omega_{b}}{\omega}$,
expression for the body acceleration takes the form
\begin{eqnarray}
G(t)=a_{0}cos(bt+\phi_{g}) ~~~~with ~~~t\geq 0
\end{eqnarray}
Also the non-dimensional form of pressure gradient and wall motion
are assumed as
\begin{eqnarray}
-\frac{\partial P(t)}{\partial z}=\bar{K}+K_{p}cos(t)
\end{eqnarray}
and
\begin{eqnarray}
 R(t)=\bar{R}[1+K_{r}sin(t+\phi_{r})]
 \end{eqnarray}
 where $\bar{R}$, $K_{r}$, $\phi_{r}$, $\bar{K}$, and $K_{p}$ stand
 for the mean radius, amplitude of arterial wall motion, the phase
 difference, constant amplitude of pressure gradient and amplitude
 of the pulsatile component gives rise to systolic and diastolic
 pressure.
 It is clear from equations (8)-(11) that when $K$ and $m$ are set to zero,
 the system of equations reduces to a classical Newtonian fluid model.\\
 Let us simplify the equations of motion by using long wavelength
 approximation $(\frac{R}{\lambda}<<1)$, where the lumen
 of the arterial radius $R$ is sufficiently smaller than the wavelength $\lambda$
 of pressure wave. Under this assumption the axial viscous transport terms are
 negligible from (8)-(11). Therefore the equation (8) reduces to simply $\frac{\partial p}{\partial
 r}=0$ and others becomes
\begin{eqnarray}
          \frac{\partial u}{\partial z}+\frac{v}{r}+\frac{\partial v}{\partial r}=0
\end{eqnarray}
\begin{eqnarray}
   \alpha^{2}(\frac{\partial u}{\partial t}+u\frac{\partial u}{\partial z}+v\frac{\partial u}{\partial
         r})=-\frac{\partial p}{\partial z}+(1+K)(\frac{\partial^{2} u}{\partial {r}^{2}}
         +\frac{1}{r}\frac{\partial u}{\partial r})+K(\frac{\partial w}{\partial
         r}+\frac{w}{r})-H^{2}u+G(t)
\end{eqnarray}
\begin{eqnarray}
   \alpha^{2} J(\frac{\partial w}{\partial t}+u\frac{\partial w}{\partial z}+v\frac{\partial w}{\partial
         r})=-K(2w+\frac{\partial u}{\partial r}-\frac{\partial v}{\partial z})+m(\frac{\partial^{2} w}{\partial r^{2}}+\frac{1}{r}\frac{\partial w}{\partial r}-\frac{w}{r^{2}})
\end{eqnarray}
\begin{eqnarray}
     \frac{\partial \theta}{\partial t}+u\frac{\partial \theta}{\partial z}+v\frac{\partial \theta}{\partial
         r}=\frac{1}{\alpha^{2}P_{r}}(\frac{\partial^{2} \theta}{\partial r^{2}}+\frac{1}{r}\frac{\partial \theta}{\partial r})+\frac{E_{c}H^{2}}{\alpha^{2}}u^{2}
\end{eqnarray}
The boundary conditions for the present problem are described as
follows:\\
Along the central line of the tube (at $r=0$) the radial velocity,
axial velocity gradient and spin of microrotation vanishes,
\begin{eqnarray}
i.e,~~ at~~ r=0,~~v=w=\frac{\partial u}{\partial r}=\frac{\partial \theta}{\partial r}=0.
\end{eqnarray}
 and at the wall
 \begin{eqnarray} i.e,~~ at~~ r=R(z,t),~~u=w=0~~~\theta=1.0~~ and~~~
v=\frac{\partial R}{\partial t}.
\end{eqnarray}
The geometry of the stenosis is described mathematically in
non-dimensional form as (cf. Fig. 1)
\begin{eqnarray}
R(z,t)&=&\bar{R}\left[1-\frac{\delta}{2R_{0}}\{1+cos(\frac{2\pi}{l_{0}}(z-d-\frac{l_{0}}{2}))\}\right](1+K_{r}sin(t+\phi_{r}));
~~~when~~ d<z<d+l_{0},\nonumber\\
      &=&\bar{R}\left[1+K_{r}sin(t+\phi_{r})\right];
      ~~~~~ elswhere
\end{eqnarray}
where $\bar{R}$ denotes the mean radius of the artery, $\delta $
the depth of the stenosis, $R_0$ is the radius of the artery at
normal pathological state.

\section{Method of Solution}
Since the governing equations are couple and non-linear in nature,
finding analytic solution is impossible. To simplify the governing
equations, we introduced a radial coordinate transformation given
by
  \[\xi=\frac{r}{R(z,t)}.\]
Using this transformation, the equations (15)-(18) can be written
in the following form
\begin{eqnarray}
          \frac{\partial u}{\partial z}-\frac{\xi}{R}\frac{\partial u}{\partial \xi}\frac{\partial R}{\partial z}+\frac{v}{\xi R}+\frac{1}{R}\frac{\partial v}{\partial
          \xi}=0,
\end{eqnarray}
\begin{eqnarray}
\frac{\partial u}{\partial t}=\frac{\frac{\partial u}{\partial
\xi}}{R}\left[\xi(\frac{\partial R}{\partial t}+u\frac{\partial
R}{\partial z})-v\right]-u\frac{\partial u}{\partial
z}+\frac{(1+K)}{\alpha^{2}R^{2}}\left[\frac{\partial^{2}
u}{\partial \xi^{2}}+\frac{1}{\xi}\frac{\partial u}{\partial
\xi}\right]+\frac{K}{\alpha^{2}R}(\frac{\partial w}{\partial
\xi}+\frac{w}{\xi})\nonumber\\+\frac{1}{\alpha^{2}}\left[G(t)-H^{2}u-\frac{\partial
p}{\partial z}\right],
\end{eqnarray}
\begin{eqnarray}
\frac{\partial w}{\partial t}=\frac{\frac{\partial w}{\partial
\xi}}{R}\left[\xi(\frac{\partial R}{\partial t}+u\frac{\partial
R}{\partial z})-v\right]-u\frac{\partial w}{\partial
z}-\frac{K}{\alpha^{2}J}(2w-\frac{\partial v}{\partial
z})-\frac{K}{\alpha^{2}JR}\left[\frac{\partial u}{\partial
\xi}+\xi\frac{\partial v}{\partial \xi}\frac{\partial R}{\partial
z}\right]\nonumber\\+\frac{m}{\alpha^{2}JR^{2}}\left[\frac{\partial^{2}w}{\partial
\xi^{2}}+\frac{1}{\xi}\frac{\partial w}{\partial
\xi}-\frac{w}{\xi^{2}}\right],
\end{eqnarray}
\begin{eqnarray}
\frac{\partial \theta}{\partial t}=\frac{\frac{\partial
\theta}{\partial \xi}}{R}\left[\xi(\frac{\partial R}{\partial
t}+u\frac{\partial R}{\partial z})-v\right]-u\frac{\partial
\theta}{\partial
z}+\frac{1}{\alpha^{2}P_{r}R^{2}}\left[\frac{\partial^{2}
\theta}{\partial \xi^{2}}+\frac{1}{\xi}\frac{\partial
\theta}{\partial \xi}\right]+\frac{E_{c}H^{2}}{\alpha^{2}}u^{2},
\end{eqnarray}
The boundary conditions (19) and (20) are also transformed to
\begin{eqnarray}
v=w=\frac{\partial u}{\partial \xi}=\frac{\partial \theta}{\partial \xi}=0~~~ at~~ \xi=0.
\end{eqnarray}
and
\begin{eqnarray}
u=w=0,~~v=\frac{\partial R}{\partial t}~~and~~\theta=1.0~~ at~~ \xi=1.
\end{eqnarray}
Let us multiply the equation (22) by $\xi R$ and integrating with
respect to $\xi$ between the limits $0$ and $\xi$, we get
\begin{eqnarray}
 v(\xi, z, t)=\xi u\frac{\partial R}{\partial
 z}-\frac{R}{\xi}\int^{\xi}_{0}\xi\frac{\partial u}{\partial z}d\xi-\frac{2}{\xi}\frac{\partial R}{\partial
 z}\int^{\xi}_{0}\xi ud\xi.
 \end{eqnarray}
 Using the boundary condition (27) at $\xi=1$ the equation (28) yeilds\\
 \[\frac{1}{R}\frac{\partial R}{\partial t}=-\frac{2}{R}\frac{\partial R}{\partial
 z}\int^{1}_{0}\xi ud\xi-\int^{1}_{0}\xi\frac{\partial u}{\partial
 z}d\xi\]
 and has the form
 \begin{eqnarray}
 \int^{1}_{0}\xi\frac{\partial u}{\partial z}d\xi=\int^{1}_{0}\xi\left[\frac{1}{R}\frac{\partial R}{\partial t}f(\xi)-\frac{2}{R}\frac{\partial R}{\partial
 z}u\right]d\xi.
 \end{eqnarray}
 Let us choose  ~$f(\xi)=-4(\xi^{2}-1)$~ such that
 $\int^{1}_{0}\xi f(\xi)d\xi=1.$\\
Now equating the integrands from both sides of (29) we have
\begin{eqnarray}
 \frac{\partial u}{\partial z}=\frac{4(\xi^{2}-1)}{R}\frac{\partial R}{\partial t}-\frac{2}{R}\frac{\partial R}{\partial
 z}u,
 \end{eqnarray}
which on substitution in (28), we obtain the radial velocity in
the following form
 \begin{eqnarray}
 v(\xi, z, t)=\xi\left[u\frac{\partial R}{\partial
 z}+(2-\xi^{2})\frac{\partial R}{\partial t}\right ].
 \end{eqnarray}
 In the next section, we present computational scheme for
 computing numerical solution of the velocity components.
\section{Computational Scheme}
 Since our Physical domain is highly complex to reduce it to a simple
 rectangular domain, we used the co-ordinate transformation, presented in the previous
 section. To solve the equations (22)-(25) and (30) using a finite difference technique,
 we subdivided the rectangular domain into a network by drawing straight lines parallel
 to the co-ordinate axes. The solutions $u(z,\xi,t)$, $v(z,\xi,t)$, $w(z,\xi,t)$ and $\theta(z,\xi,t)$
 at any mesh point in a computational domain are denoted by $u^{k}_{i,j}$, $v^{k}_{i,j}$,
 $w^{k}_{i,j}$ and $\theta^{k}_{i,j}$ in which\\
 $z_{i}=i\Delta z$~;~~ $i=0, 1, \cdots , M$;\\
 $\xi_{j}=j\Delta \xi$~;~~ $j=0, 1, \cdots ,N$;\\
 $t_{k}=(k-1)\Delta t$~;~~ $k= 1, 2, 3 ......$\\
where $M$ and $N$ are the maximum number of mesh points in the $z$
and $\xi$-directions respectively. Let $V^{k}_{i,j}$ stand for the
variables $u$, $v$, $w$ and $\theta$, then the spatial and time
derivatives are replaced by the central difference as\\
$(\frac{\partial V}{\partial \xi})^{k}_{i,j}=\frac{V^{k}_{i,
j+1}-V^{k}_{i, j-1}}{2\Delta
 \xi}$,~ $(\frac{\partial V}{\partial z})^{k}_{i,j}=\frac{V^{k}_{i+1,j}-V^{k}_{i-1,j}}{2\Delta z}$,~
$(\frac{\partial^{2} V}{\partial
\xi^{2}})^{k}_{i,j}=\frac{V^{k}_{i,j+1}-2V^{k}_{i,j}+V^{k}_{i,j-1}}{(\Delta
\xi)^{2}}$,~ $(\frac{\partial V}{\partial
t})^{k}_{i,j}=\frac{V^{k+1}_{i, j}-V^{k}_{i, j}}{\Delta t}$ and so on.\\
Using these finite differences, the discretized equations for
(23)-(25) and (30) becomes
\begin{eqnarray}
u^{k+1}_{i, j}=u^{k}_{i, j}+\Delta t
[\frac{u^{k}_{i,j+1}-u^{k}_{i, j-1}}{2\Delta \xi
R^{k}_{i}}\left(\xi_{j}\left(\frac{\partial R}{\partial
t}^{k}_{i}+u^{k}_{i, j}\frac{\partial R}{\partial
z}^{k}_{i}\right)-v^{k}_{i, j}\right)-u^{k}_{i,
j}\frac{u^{k}_{i+1, j}-u^{k}_{i-1, j}}{2\Delta
z}+\nonumber\\
\frac{1+K}{\alpha^{2}(R^{k}_{i})^{2}}\left(\frac{u^{k}_{i,
j+1}-2u^{k}_{i,j}+u^{k}_{i, j-1}}{(\Delta
\xi)^{2}}+\frac{1}{\xi_{j}}\frac{u^{k}_{i, j+1}-u^{k}_{i,
j-1}}{2\Delta
\xi}\right)+\frac{K}{\alpha^{2}R^{k}_{i}}(\frac{w^{k}_{i,
j+1}-w^{k}_{i, j-1}}{2\Delta \xi}+\nonumber\\
\frac{w^{k}_{i,j}}{\xi_{j}})+\frac{1}{\alpha^{2}}\left(G(t)-\frac{\partial
p}{\partial z}^{k}_{i}-H^{2}u^k_{i,j}\right)],
\end{eqnarray}
\begin{eqnarray}
w^{k+1}_{i, j}=w^{k}_{i, j}+\Delta t[\frac{w^{k}_{i,j+1}-w^{k}_{i,
j-1}}{2\Delta \xi R^{k}_{i}}\left(\xi_{j}\left(\frac{\partial R}{\partial
t}^{k}_{i}+u^{k}_{i, j}\frac{\partial R}{\partial
z}^{k}_{i}\right)-v^{k}_{i, j}\right)-u^{k}_{i, j}\frac{w^{k}_{i+1,
j}-w^{k}_{i-1, j}}{2\Delta
z}-\nonumber\\ \frac{K}{\alpha^{2}J}\left(2w^{k}_{i,j}-\frac{v^{k}_{i+1,j}-v^{k}_{i-1,j}}{2\Delta z}\right)-\frac{K}{\alpha^{2}R^{k}_{i}J}\left(\frac{u^{k}_{i,
j+1}-u^{k}_{i, j-1}}{2\Delta \xi}+\xi_{j}\frac{\partial
R}{\partial z}^{k}_{i}\frac{v^{k}_{i, j+1}-v^{k}_{i, j-1}}{2\Delta
\xi}\right)+\nonumber\\ \frac{m}{\alpha^{2}(R^{k}_{i})^{2}J}\left(\frac{w^{k}_{i,
j+1}-2w^{k}_{i,j}+w^{k}_{i, j-1}}{(\Delta \xi)^{2}}+
\frac{1}{\xi_{j}}\frac{w^{k}_{i, j+1}-w^{k}_{i, j-1}}{2\Delta
\xi}-\frac{w^{k}_{i, j}}{\xi^{2}_{j}}\right)],
\end{eqnarray}
\begin{eqnarray}
\theta^{k+1}_{i, j}=\theta^{k}_{i, j}+\Delta
t[\frac{\theta^{k}_{i,j+1}-\theta^{k}_{i, j-1}}{2\Delta \xi
R^{k}_{i}}\left(\xi_{j}\left(\frac{\partial R}{\partial t}^{k}_{i}+u^{k}_{i,
j}\frac{\partial R}{\partial z}^{k}_{i}\right)-v^{k}_{i, j}\right)-u^{k}_{i,
j}\frac{\theta^{k}_{i+1, j}-\theta^{k}_{i-1, j}}{2\Delta
z}\nonumber\\+\frac{1}{\alpha^{2}P_{r}(R^{k}_{i})^{2}}\left(\frac{\theta^{k}_{i,j+1}-2\theta^{k}_{i,j}+\theta^{k}_{i,
j-1}}{(\Delta
\xi)^{2}}+\frac{1}{\xi_{j}}\frac{\theta^{k}_{i,j+1}-\theta^{k}_{i,j-1}}{2\Delta
\xi}\right)+\frac{E_{c}H^{2}}{\alpha^{2}}(u^k_{i,j})^{2}],
\end{eqnarray}
\begin{eqnarray}
v^{k+1}_{i, j}=\xi_{j}\left[u^{k}_{i, j}\frac{\partial R}{\partial
z}^{k}_{i}+(2-\xi^{2}_{j})\frac{\partial R}{\partial t}^{k}_{i}\right]
\end{eqnarray}
where $\left(\frac{\partial R}{\partial z}\right)^{k}_{i}$ and
$\left(\frac{\partial
R}{\partial t}\right)^{k}_{i}$ are computed at the arterial wall.\\
After having determined the values of the axial velocity, radial
velocity, microrotation component and temperature one can easily
obtain the fluid acceleration $F$, the volumetric flow rate $Q$,
flow resistance $\lambda$, wall shear stress $\tau_{w}$ and the
Nusselt number $Nu$ from the following relations,\\
$F^{k+1}_{i, j}=\frac{u^{k+1}_{i,j}-u^{k}_{i, j}}{\Delta t}
+u^{k}_{i,j}\frac{u^{k}_{i+1,j}-u^{k}_{i-1,j}}{2\Delta z}- \frac{\xi_{j}}{R^{k}_{i}}\left(\frac{\partial R}{\partial t}^{k}_{i}+u^{k}_{i,
j}\frac{\partial R}{\partial z}^{k}_{i}\right)\frac{u^{k}{i,j+1}-u^{k}_{i,
j-1}}{2\Delta \xi}$\\
$Q_{i}=2\pi(R^{k}_{i})^{2}\int^{1}_{0}\xi_{j}u^{k}_{i,
j}d\xi_{j}$\\
 $\lambda^{k}_{i}=\frac{|L\frac{\partial p}{\partial
z}^{k}_{i}|}{Q^{k}_{i}},$~\\
$Nu_{i}=\frac{1}{R^{2}_{0}R^{2}_{i}}\left[\frac{\theta^{k}_{i,N}-\theta^{k}_{i,N-1}}{\Delta
\xi}\right]$

\section{Computational Results and Discussion}
The objective of the present study has been to investigate the
flow characteristics of blood through stenosed artery under the
action of a periodic body acceleration as well as an external
magnetic field on the axial velocity by taking into account the
effect of microrotation of blood cells. A specific numerical
illustration has been taken to examining the applicability of the
physiological data available in the existing scientific
literatures [3,4,19-21,28,29]. The computational work has been
carried out by using the
following data:\\
$L=5.0,~ d=2.0,~l_{0}=1.0,~\delta=0.1, 0.25, 0.50; a_{0}=0.0, 1.0,
2.0, 3.0; ~ b=1.0,~ \phi_{g}=0.0,~ \phi_{r}=0.0,~\\
~ R_{0}=1.0,~\bar{R}=1.0,~K_{p}=1.46,~ \bar{K}=7.30,~
K_{r}=0.05,~f_{p}=1.2,~K=0.0, 0.1, 0.2, 0.3;~ \\
J=0.1,~m=0.1, 0.01, 0.001;~\alpha=3.0,~ H=0.0, 1.0, 2.0, 3.0, 4.0; ~E_{c}=0.0002,$\\
and for a human body temperature, $T=310K$ the value of $P_{r}=21$
is considered for blood. For the sake of comparison, we have also
examined the cases where $P_{r}=7, 14$. We used $\Delta \xi =
0.025$, $\Delta z= 0.05$, $\Delta t=0.001$ through out the
computation. We observed that further reduction in $\Delta \xi$
and $\Delta z$ does not bring about any significant change, which
leads to the stability of the numerical methods.\\

Figs. 2(a)-(c) give the variation of axial velocity $u$ along the
radial direction $\xi$. In Fig. 2(a) dotted lines indicate the
axial velocity distribution for Newtonian fluid model and solid
lines represent that of the micropolar fluid model. From this
figure, we observe that the axial velocity decreases with the
increase of the Hartmann number $H$. The axial velocity occurs
maximum at the central line of the artery in all four cases. One
can note from Fig 2(a) that for a sufficient strength of magnetic
field, the profiles for both the Newtonian and non-Newtonian
models are same. The axial velocity is found to be greater in the
case of Newtonian model than that of the micropolar fluid model.
This happens due to the presence of microrotation of blood cells,
creates an additional viscosity called rotational viscosity, which
in turn diminishes axial velocity. Fig. 2(b) depicts the variation
of axial velocity for different values of $K$, the ratio of
viscosity $\frac{k}{\mu}$. When the values of $K$ increases, the
axial velocity decreases. The results, thus obtained due to the
increase of the microrotations of blood cells in the lumen of the
artery. The axial velocity is minimum in the vicinity of the
arterial walls, which gives rise to conclude that the microtation
of blood cells occur in the plasma layer. The variation of axial
velocity for different amplitudes of body acceleration is shown in
Fig. 2(c). It reveals that the axial velocity increases as the
amplitude of body acceleration increases. Thus in the presence of
vibration environmental system may produce the increase in blood velocity.\\

Figs. 3(a)-(c) illustrate the distribution of micro-component of
velocity $w$ along the radial direction of the artery. Fig. 3(a)
shows that the microrotation component increases with the increase
of the ratio of the viscosity $K$, which in turn produces decrease
of axial velocity. It has been shown in Fig. 3(b) that the
microcomponent of velocity $w$ increases with the decrease in
micropolar material parameter $m$. The result lies in the fact
that, when m decreases, the fluid viscosity $\mu$ increases,
thereby the microcomponent of velocity increases. One can note
that the effect of microrotation observed significantly near the
arterial wall. The variation of microrotation component is found
to increase with the increase of the blockage of an artery
presented in Fig. 3(c). In all these cases $\delta =0.1,~0.25,$
and 0.5 have been considered for computation and which are
respectively represent the 10\%, 25\% and 50\% of the blockage of
the arterial segment.\\

In Figs. 4(a) and (b) we present the variation of dimensionless
temperature distribution at the throat of the stenosis and the
rate of heat transfer at the arterial wall along with the axial
direction. It is observed from Fig. 4(a) that the temperature
decreases with the increase of the Prandtl number $P_{r}$.
However, the temperature increases with the increase of Hartmann
number $H$. It is interesting to mention here that the temperature
rapidly increases with the increase of $H$ up to 4.0 (T) and
beyond which, no significant change is noticed. The rate of heat
transfer $(Nu)$ increases with the increase of prandtl number
$P_{r}$, shown in Fig. 4(b). The rate of heat transfer is maximum
at the throat of the stenosis and is minimum at the down stream of
the stenosis than that of upstream. Thus the temperature can be
increased whenever necessary by the application of the certain
magnetic field strength. It is well known that the objective of
the hyperthermia in cancer therapy is to raise the temperature
above a therapeutic value $42^{0}$ C, while
maintaining the surrounding temperature at sublethal temperature \cite{R38}.\\

The important characteristics of blood flow are the wall shear
stress $\tau_{w}$ and the volumetric flow rate $Q$. Figs. 5(a)-(d)
give the variation of wall shear stress $\tau_{w}$ with axial
distance $z$ and Time $T=t/2\pi f_{p}$ for different height of the
stenosis $\delta$, different values of Hartmann number $H$,
amplitude of body acceleration $a_{0}$ and the material parameter
$K$. In Fig. 5(a), we observe that the wall shear stress increases
as the height of the stenosis increases. It is interesting to note
from this figure that the wall shear stress significantly changes
at the down stream of the stenosis, whereas, at the upstream of
the stenosis no change is occur. Fig. 5(b) presents the variation
of wall shear stress $\tau_{w}$ with time $T$ for different values
of the Hartmann number $H$. It shows that the wall shear stress
varies periodically with time $T$ and the peak value of each
oscillation decreases with the increase of $H$. Fig. 5(c) shows
that the amplitude of the oscillation of wall shear stress
increases with the increase of the amplitude of the body
acceleration. However, there is no significant change on the wall
shear stress for the variation of material parameter $K$ (cf. Fig.
5(d)). This is due to the fact that the microelements are unable
to rotate close to the arterial wall. Although the peak value of
oscillation increases with the increase of $K$, but the change is insignificant.\\

The variation of volumetric flow rate $Q$ with time $T$ presented
through the Figs. 6(a)-(c). We observe that the flow rate $Q$
decreases with the increase of Hartmann number $H$ as well as the
increase of the material parameter $K$, while the trend is
reversed in the case of the amplitude of body acceleration $G(t)$.
In all these three cases, the flow rate $Q$ varies periodically
with time $T$. It thus turns out that under the action of a
magnetic field, volume of blood flow can be controlled during
surgeries. The volumetric flow has an enhancing effect with its
peak value of oscillation on the vibration environmental system.\\

Figs. 7(a) and (b) depict the variation of fluid acceleration $F$
along the radial coordinate $\xi$ for different values of the
amplitude of body acceleration and different magnetic field
intensity given by $H$. For each of these cases, the fluid
acceleration $F$ is maximum at the central line and then it is
gradually diminishes and ultimately vanishes at the wall of the
artery. It is also observe from Fig. 7(a) that the magnitude of
the fluid acceleration $F$ increases as the amplitude of body
acceleration increases. Thus in the presence of an environmental
system fluid acceleration enhances. However, the magnitude of the
fluid acceleration $F$ decreases with the increase of the Hartmann
number $H$, shown in Fig. 7(b). One can mention here that the rate
of increase or decrease of fluid acceleration $F$ significantly
affected by the
application of magnetic field more than in the presence of body acceleration.\\

The variation of flow resistance $\lambda$ along with the height
 of the stenosis $\delta$ has been shown in Figs. 8(a)-(c). It has
been observed from Fig. 8(a) that the resistance to flow decreases
with an increase in amplitude $a_{0}$ of the body acceleration. On
the other hand flow resistance increases with the increase of the
values of $H$ and $K$, shown in Figs. 8(b) and (c). It is worth
mentioning that in all these figures, the flow resistance
$\lambda$ gradually increases as the depth of the stenosis
increases. It is learn from these figures that the flow resistance
remains fixed upto a certain height of the stenosis that is upto
25\% of the blockage of an artery (called mild
stenosis) beyond which it increases significantly.\\

\section {Concluding Remarks}
The present mathematical analysis is motivated towards the
modelling of pulsatile blood flow through a stenosed artery under
the combined influence of an external magnetic field and periodic
body acceleration. The mathematical model is developed towards the
application in blood flow in order to consider the effect of
micro-rotation of micro-particles suspended in plasma. When a
magnetic field is applied to the blood flow, the erythrocytes
orient their disk plane parallel to the axis of magnetic field,
creates an additional viscosity, which in turn produces decrease
of blood velocity. This interesting observation can be seen in
Fig. 2(b). The temperature of blood can be increased by the
application of magnetic field. It is also concluded that after
certain magnetic field strength, the temperature is constant. This
important result can be useful in the therapeutic treatment of
patient \cite{R38}. The present study also enables to have an
estimate of the effects of microcomponent of velocity under the
influence of periodic body acceleration. However, the
consideration of pulsatile flow on the cases of the severe
stenosis are the further
scope of the study. \\

 \noindent {\bf Acknowledgement:} {\it The authors are thankful to both the
reviewers for their kind words of appreciation and nice comments
on the applicability of the results presented. The original
manuscript has been revised on the basis of the reviewers'
suggestions. One of the authors G. C. Shit is grateful to DST, New
Delhi for awarding BOYSCAST Fellowship (Ref. No. SR/BY/M-01/09) to
him.}

\newpage
\textheight 22.0cm \pagebreak
\begin{minipage}{1.0\textwidth}
   \begin{center}
      \includegraphics[width=3.8in,height=2.5in ]{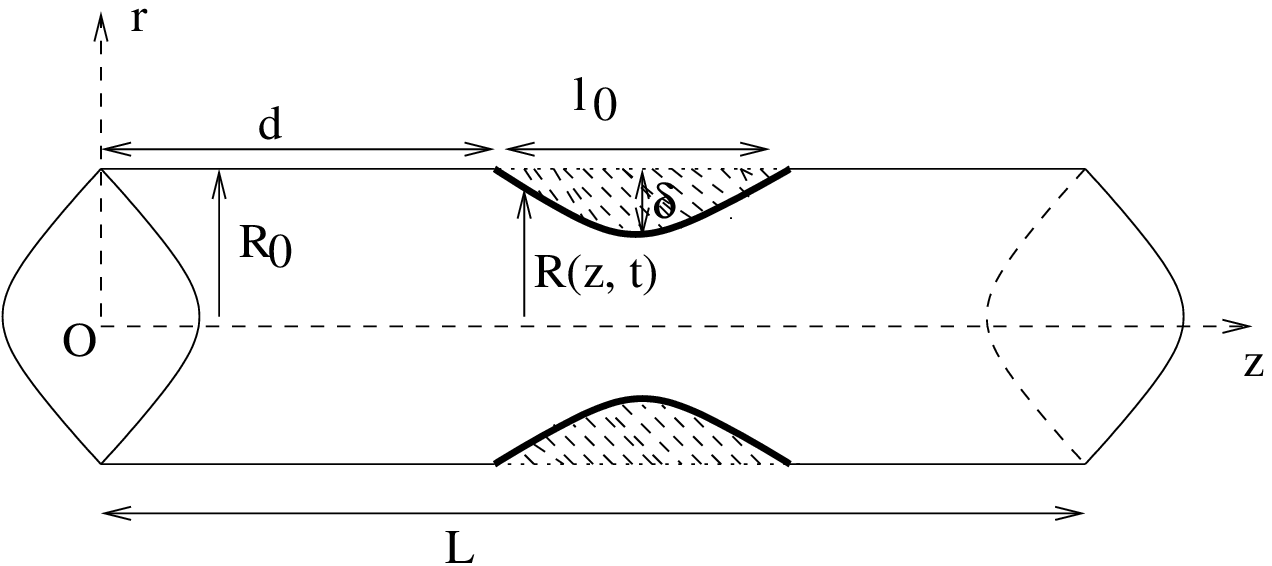}\\
Fig. 1 ~~A physical sketch of the model artery with stenosis. \\
\end{center}
\end{minipage}\vspace*{.5cm}\\
\newpage

\begin{minipage}{1.0\textwidth}
  \begin{center}

    (a)\includegraphics[width=3.0in,height=2.0in ]{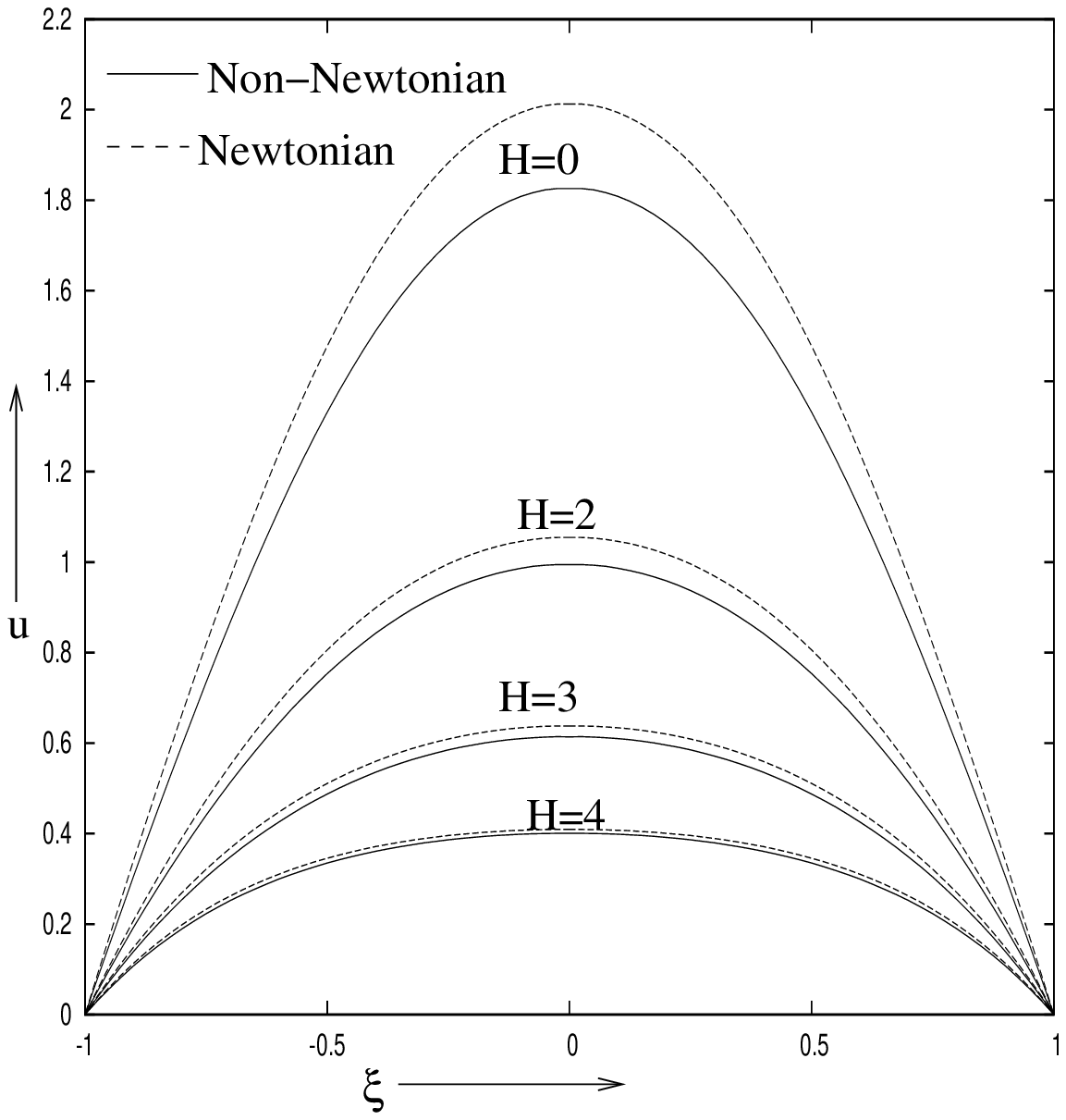}~~~
    (b)\includegraphics[width=3.0in,height=2.0in ]{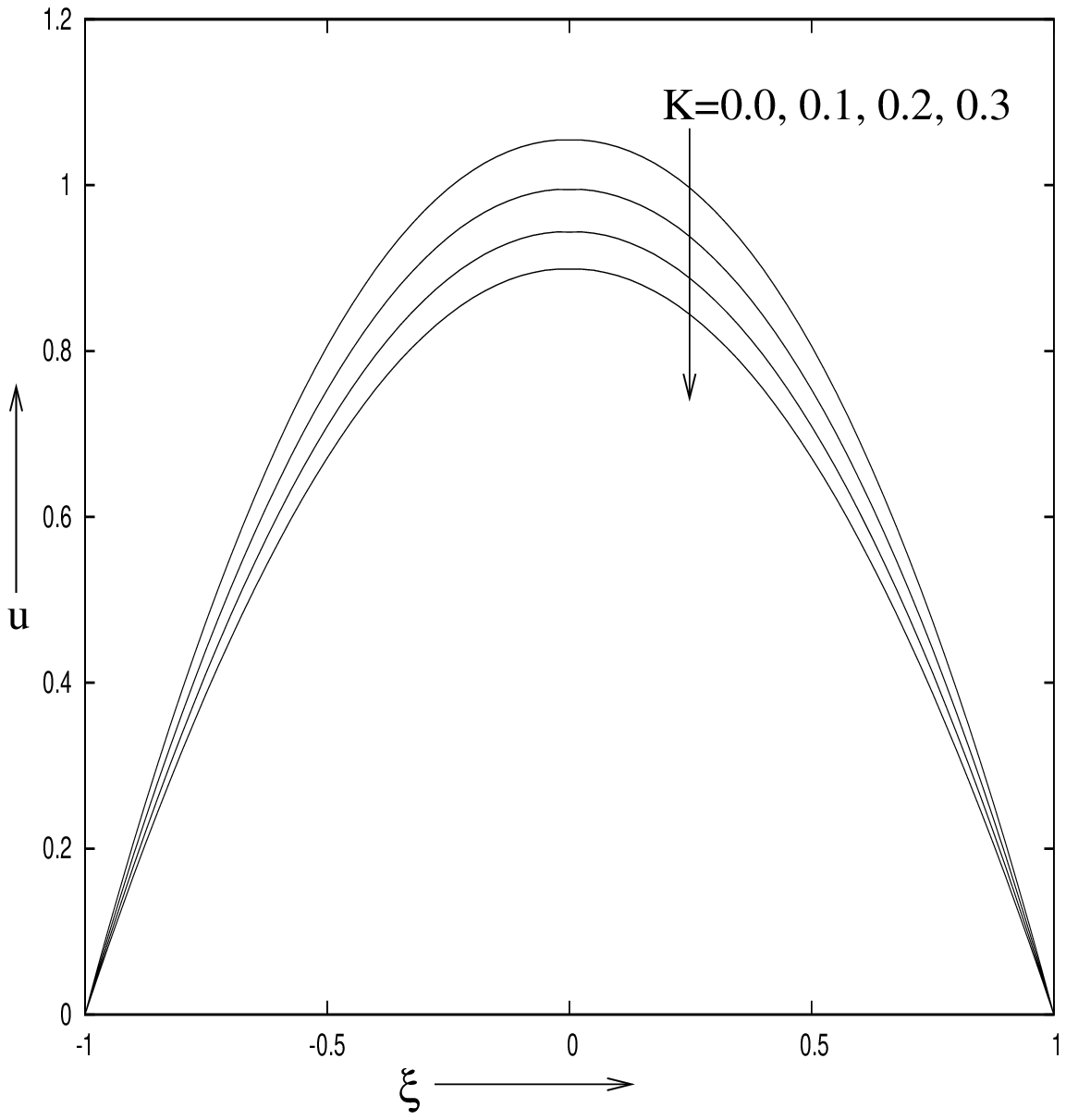}\\
    (c)\includegraphics[width=3.0in,height=2.0in ]{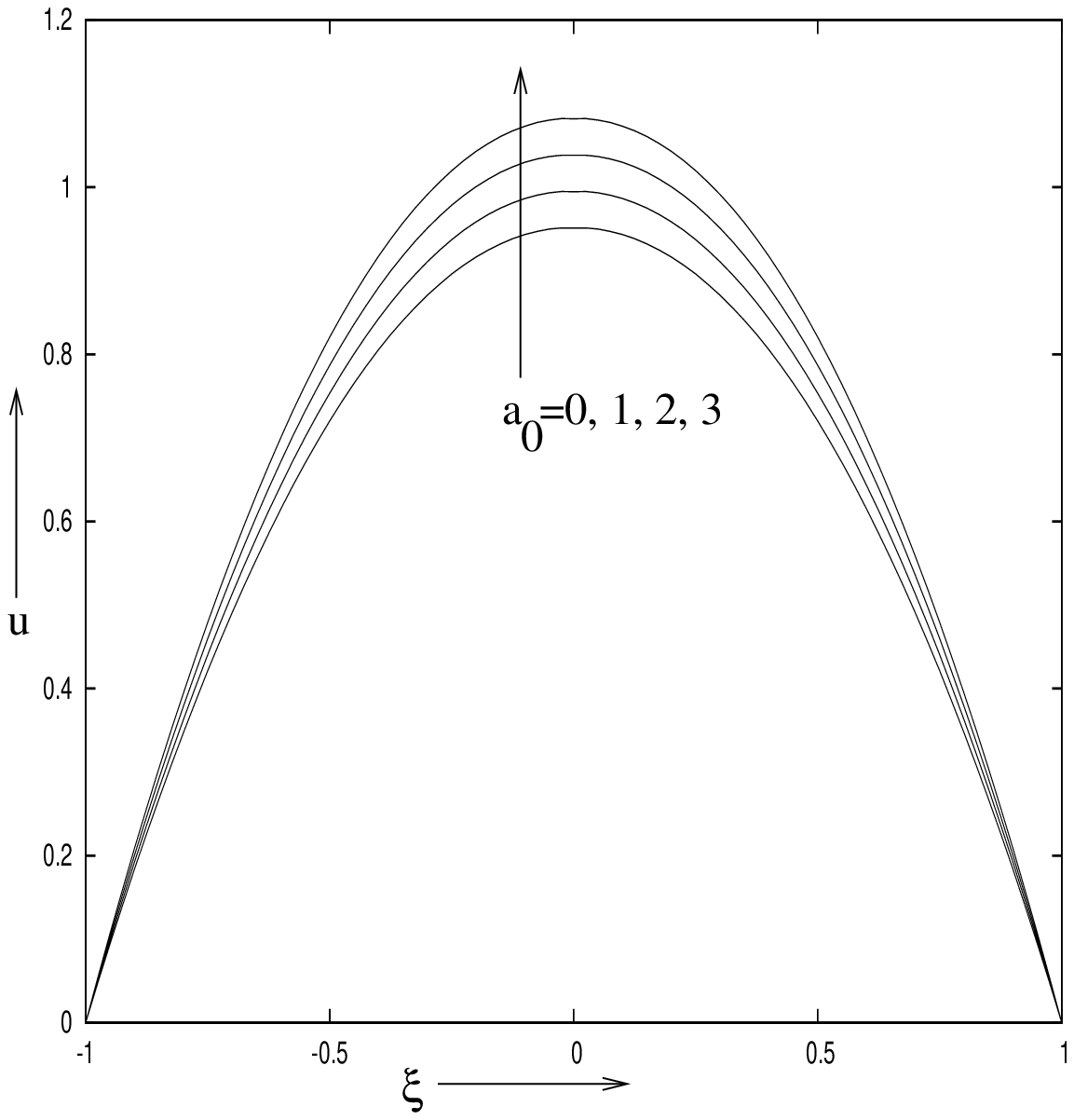}~~~\\
Fig. 2 Variation of axial velocity $u$ at the throat of the
stenosis for (a) different values of the Hartmann number $H$, when
$K=0.1$, $a_{0}=1.0$; (b) different values of $K$ when $H=2.0$,
$a_{0}=1.0$; (c) different amplitude of body
acceleration $a_{0}$ when $H=2.0$ and $K=0.1$.\\
\end{center}
\end{minipage}\vspace*{.5cm}\\
\newpage

\begin{minipage}{1.0\textwidth}
\begin{center}
     (a) \includegraphics[width=3.0in,height=2.0in]{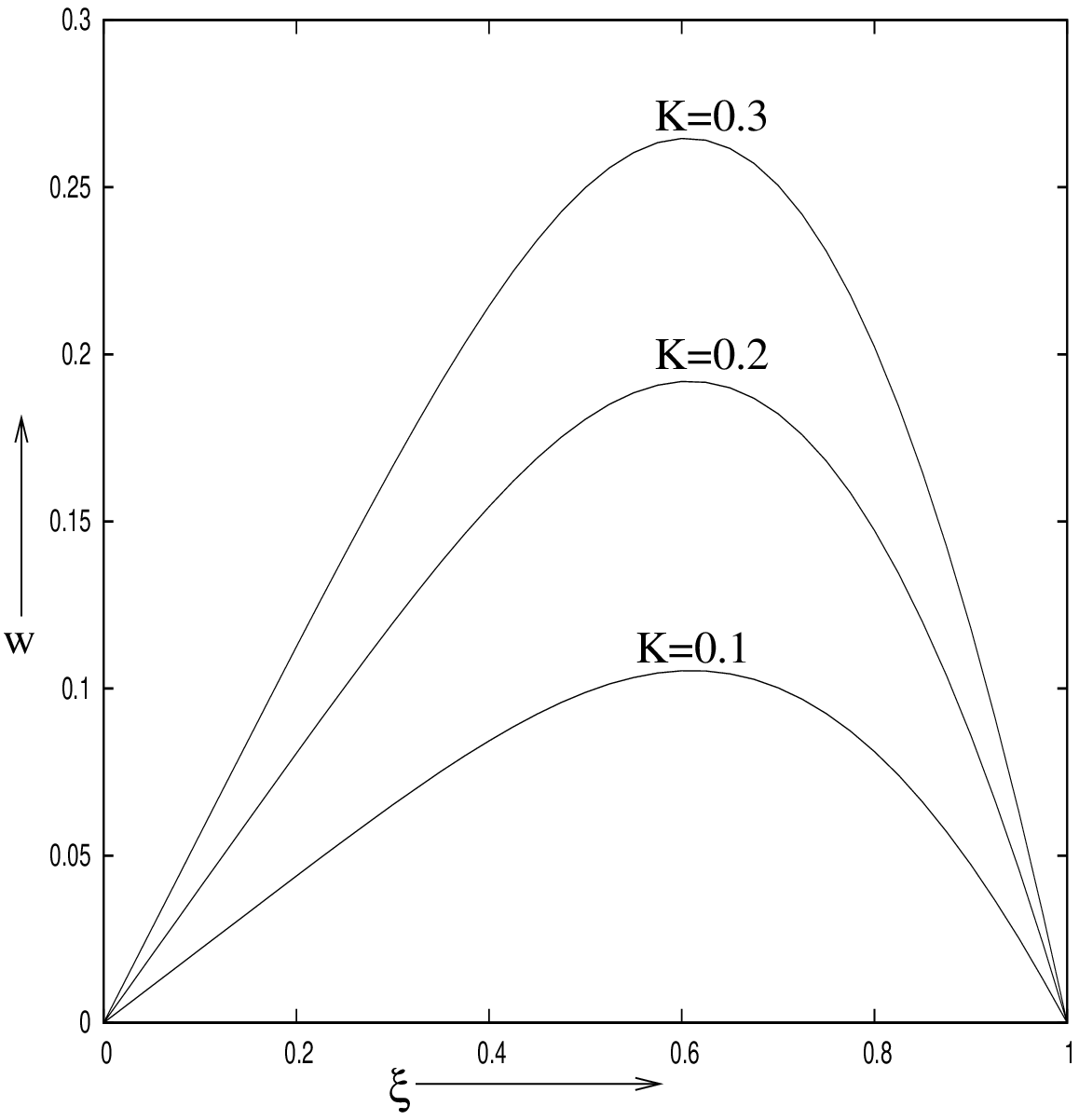}~~~
     (b) \includegraphics[width=3.0in,height=2.0in ]{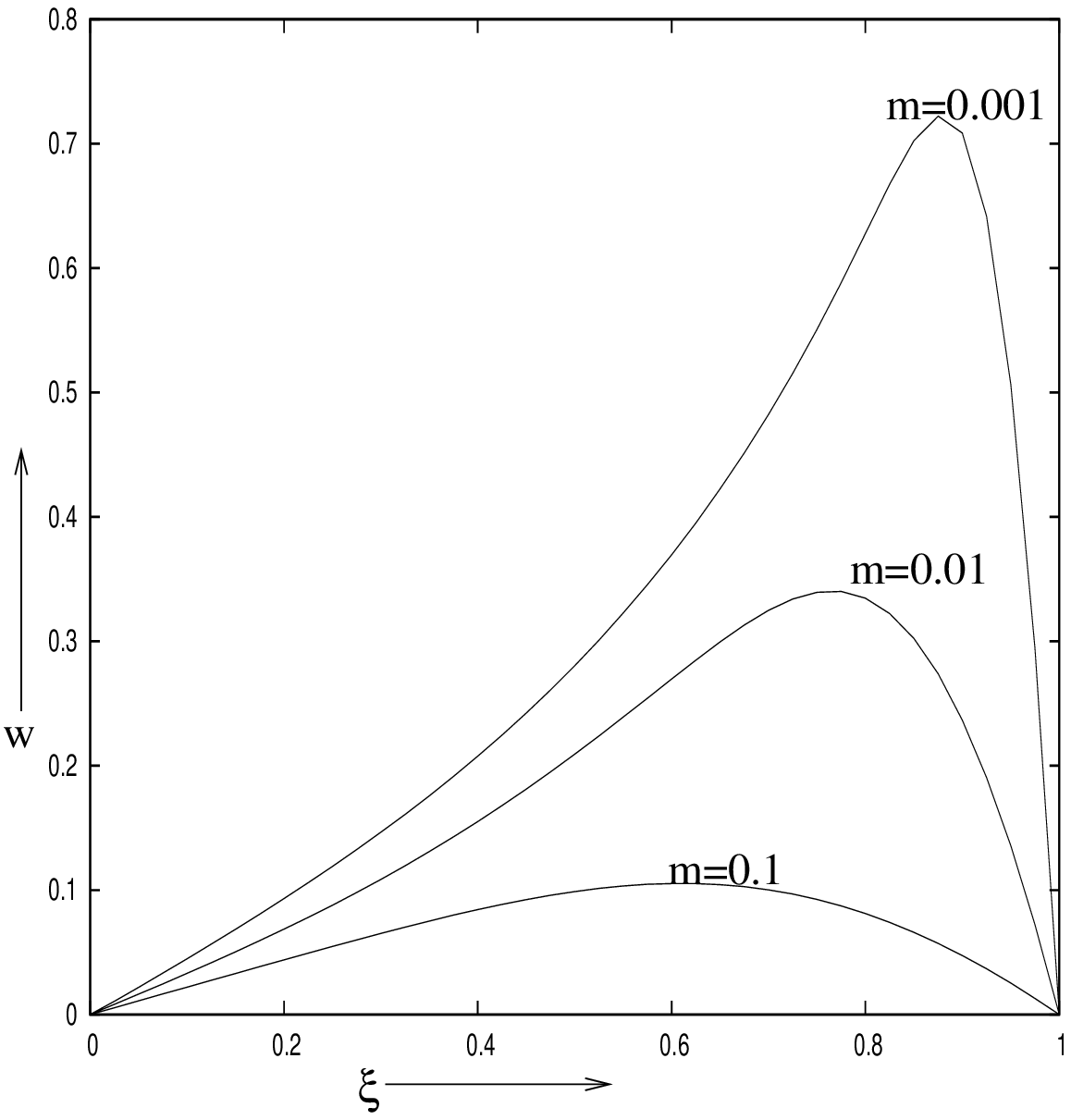}\\
     (c) \includegraphics[width=3.0in,height=2.0in ]{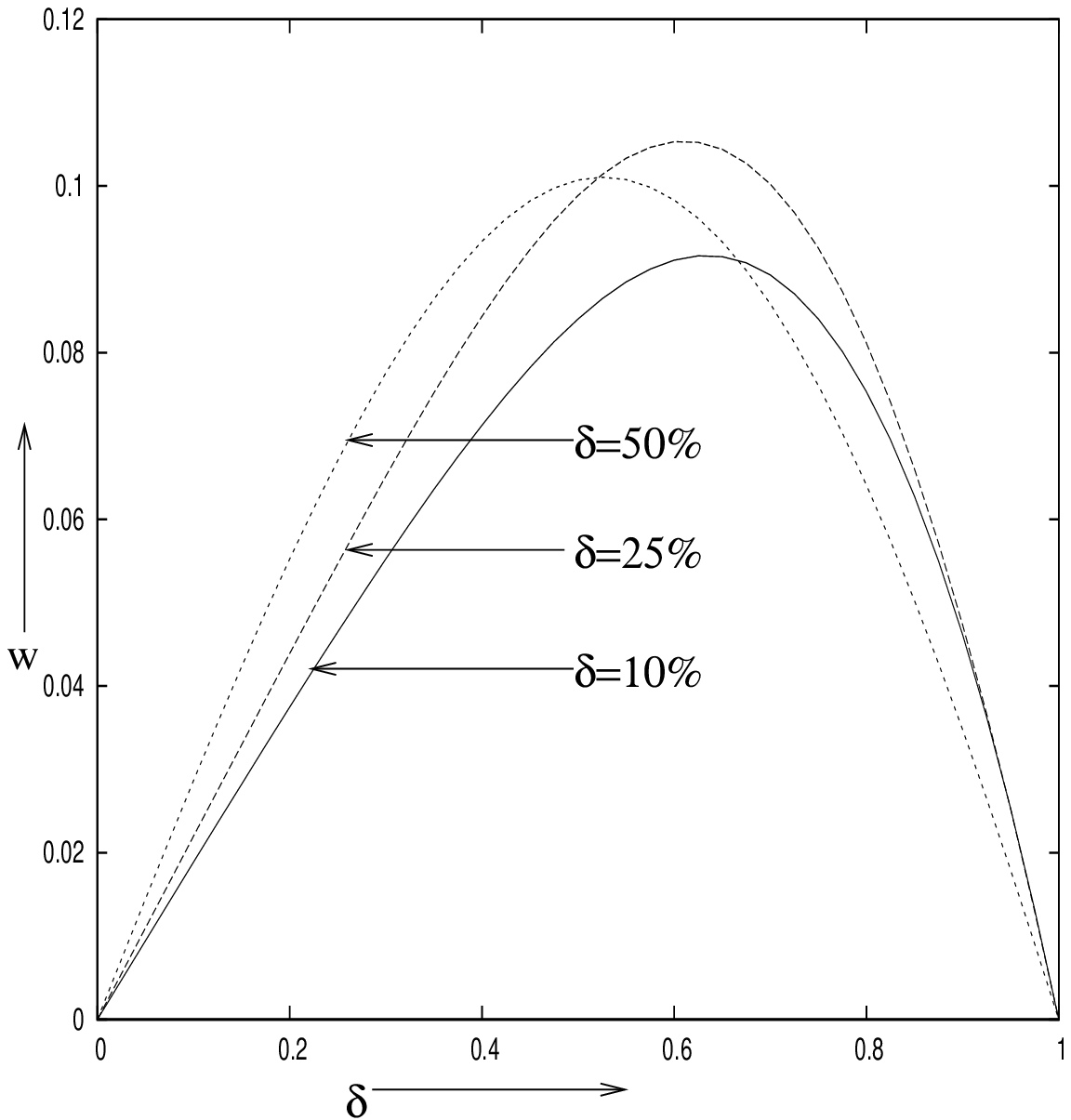}\\
Fig. 3 Variation of microrotation component $w$ at the throat of
the stenosis (a) for different values of material parameter $K$
with $m=0.1$, $\delta=0.25$; (b) for different values of $m$, with
$K=0.1$, $\delta=0.25$; (c) for different values of $\delta$, with
$K=0.1$, $m=0.1$.
\end{center}
\end{minipage}\vspace*{.5cm}\\
\newpage

\begin{minipage}{1.0\textwidth}
  \begin{center}
   (a)\includegraphics[width=3.0in,height=2.0in ]{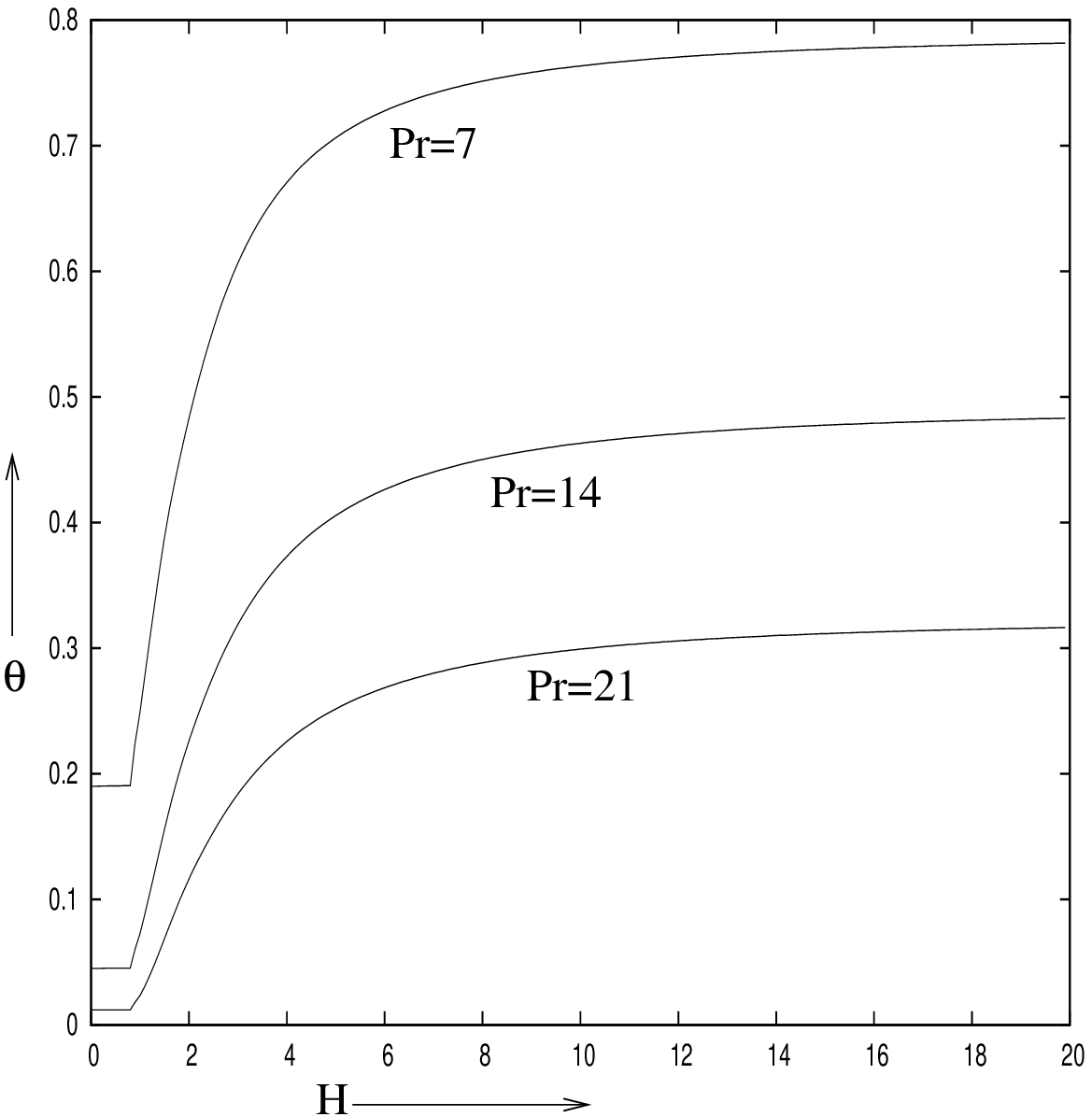}~~~
   (b) \includegraphics[width=3.0in,height=2.0in ]{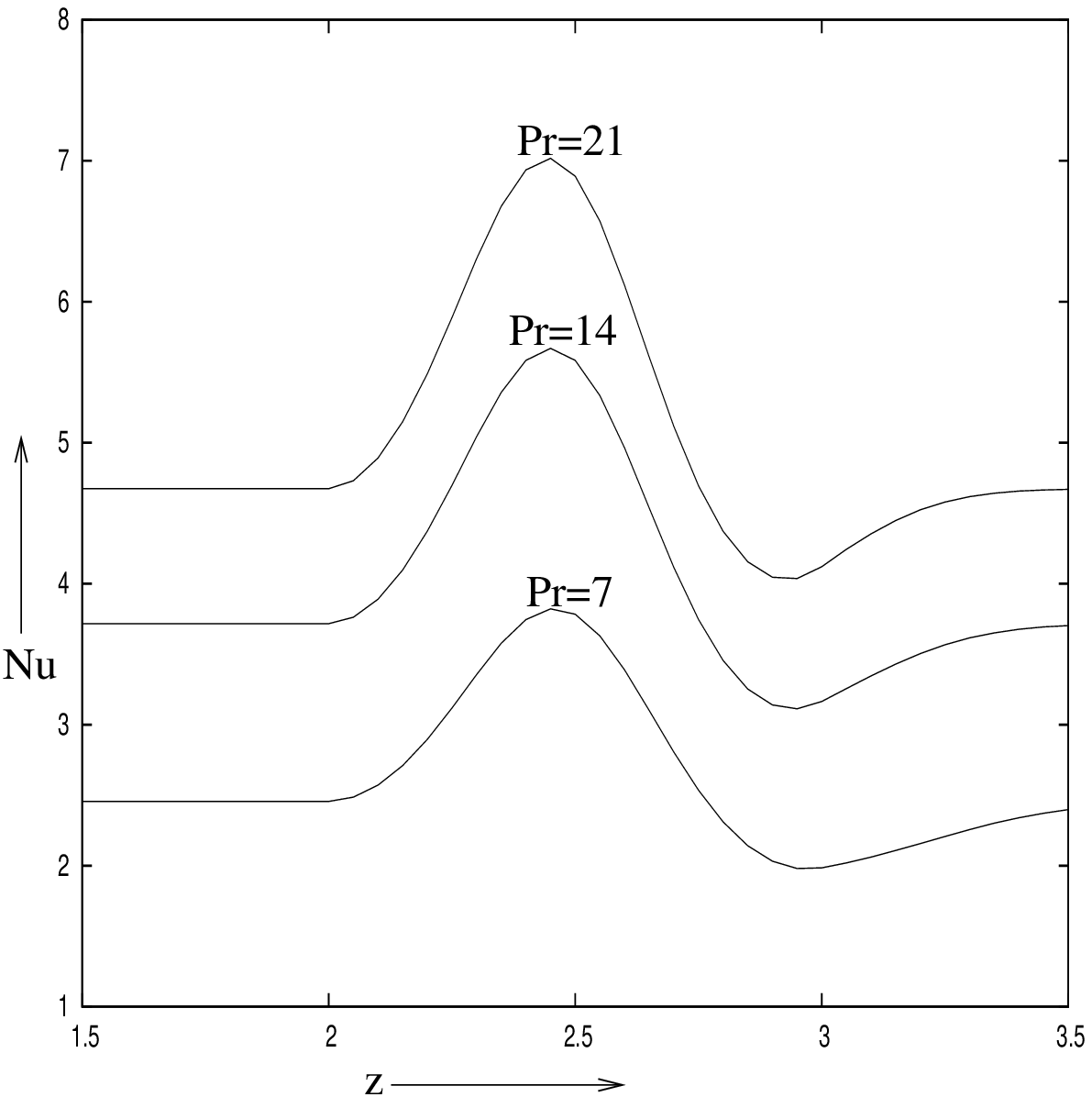}\\
Fig. 4 (a) Distribution of temperature $\theta$ for different
values of $H$ and $Pr$. (b) Variation of Nusselt number $Nu$ along
the axis of the artery for different values of Prandlt number
$Pr$. With $a_{0}=1.0$,~$b=1.0$,~$\delta=0.25$,~ $\phi_{g}=0.0,~
\phi_{r}=0.0, K=0.1,~ m=0.1,~\alpha=3.0,~ H=2.0,$
\end{center}
\end{minipage}\vspace*{.5cm}\\
\newpage

\begin{minipage}{1.0\textwidth}
   \begin{center}
(a)\includegraphics[width=3.0in,height=2.0in ]{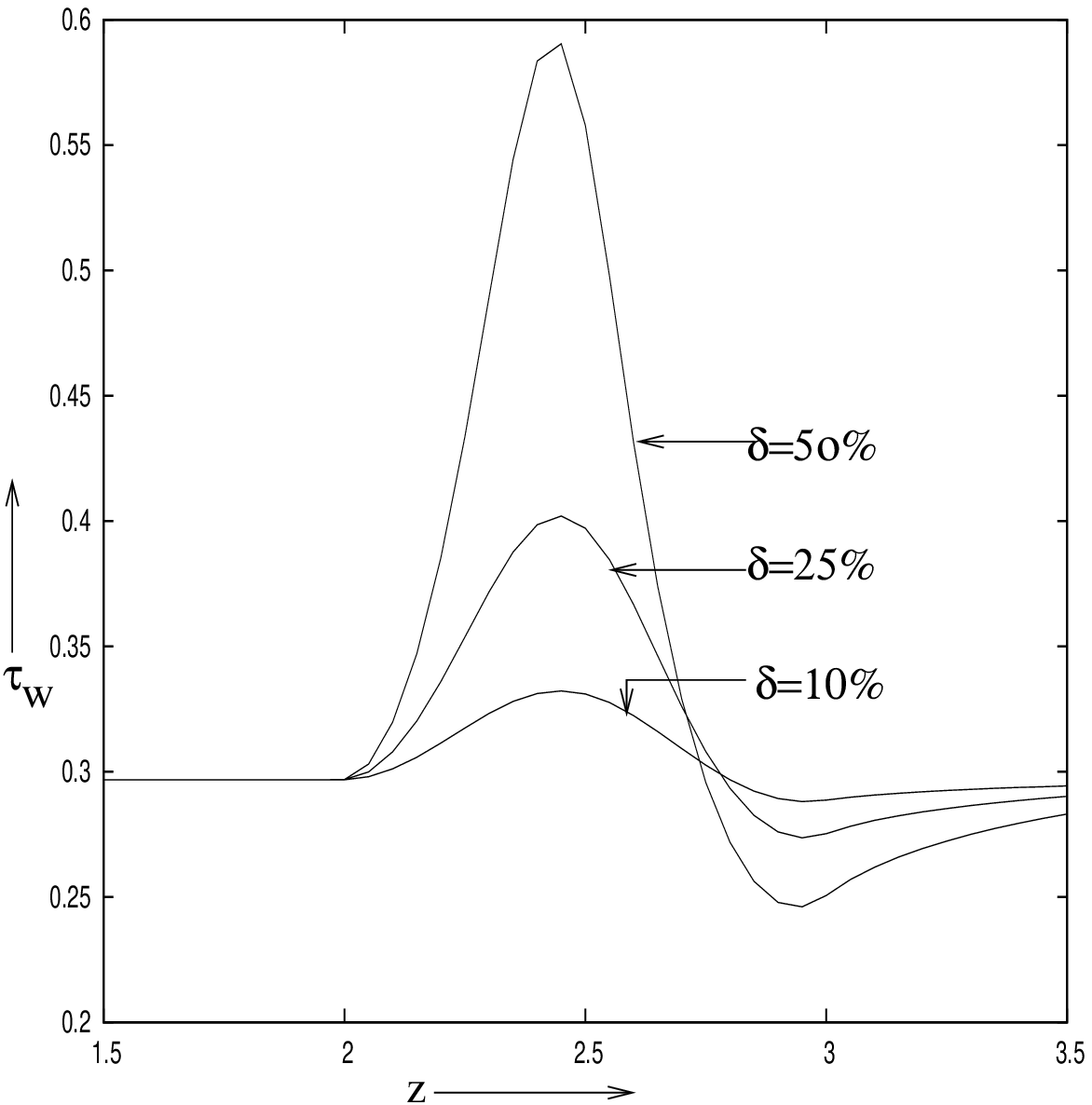}~~~
(b)\includegraphics[width=3.0in,height=2.0in ]{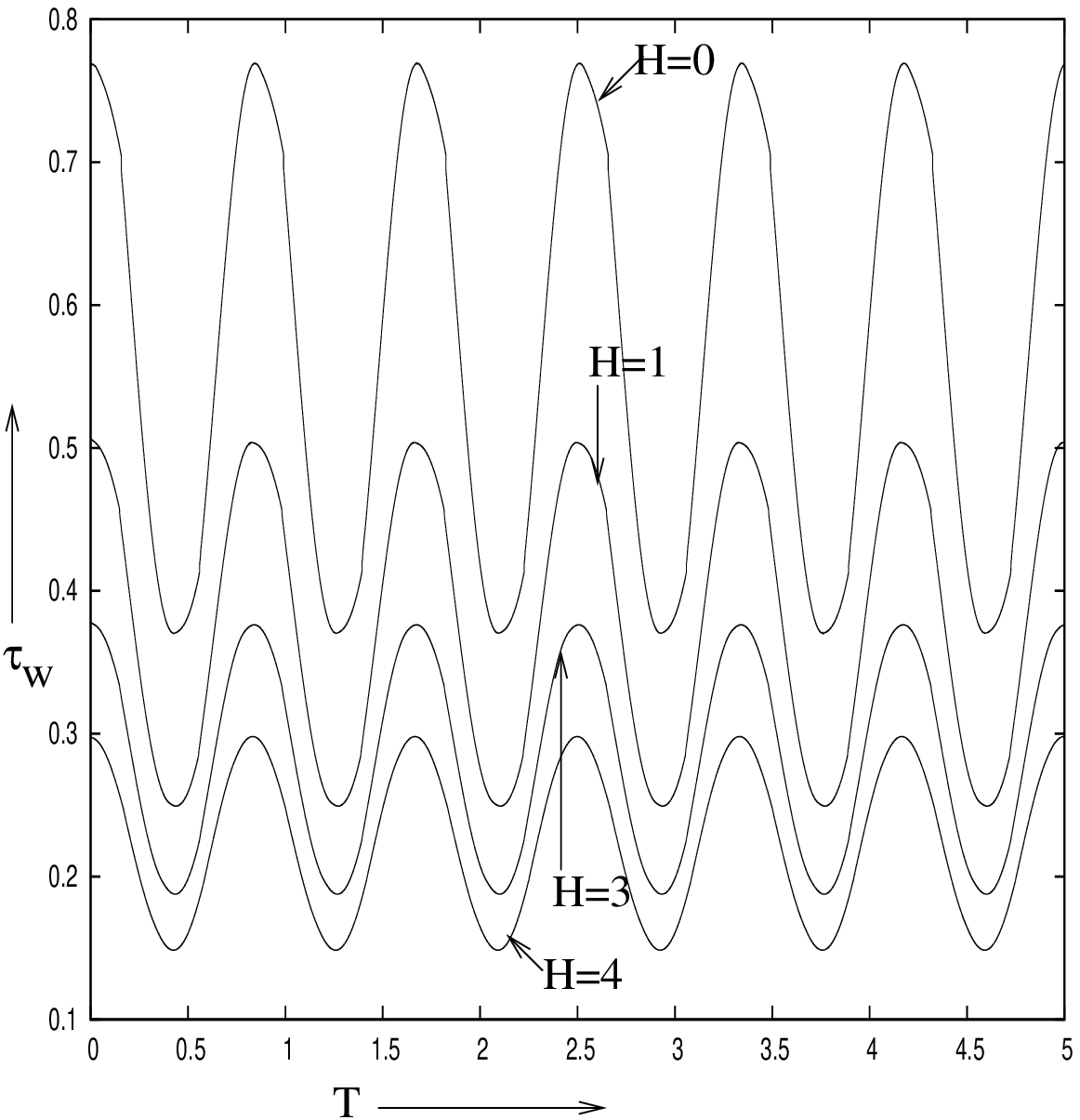}\\
(c)\includegraphics[width=3.0in,height=2.0in ]{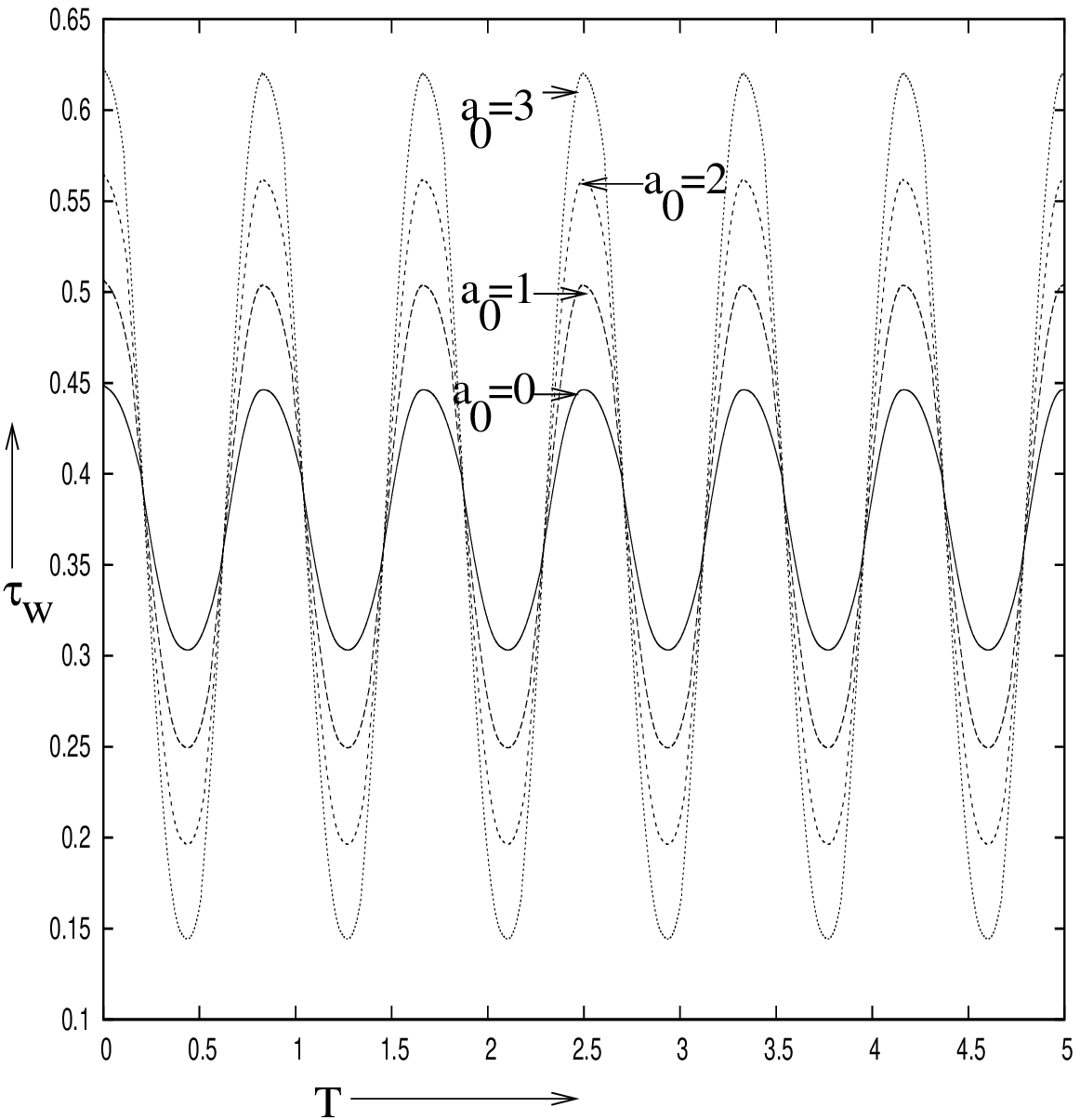}~~~
(d)\includegraphics[width=3.0in,height=2.0in ]{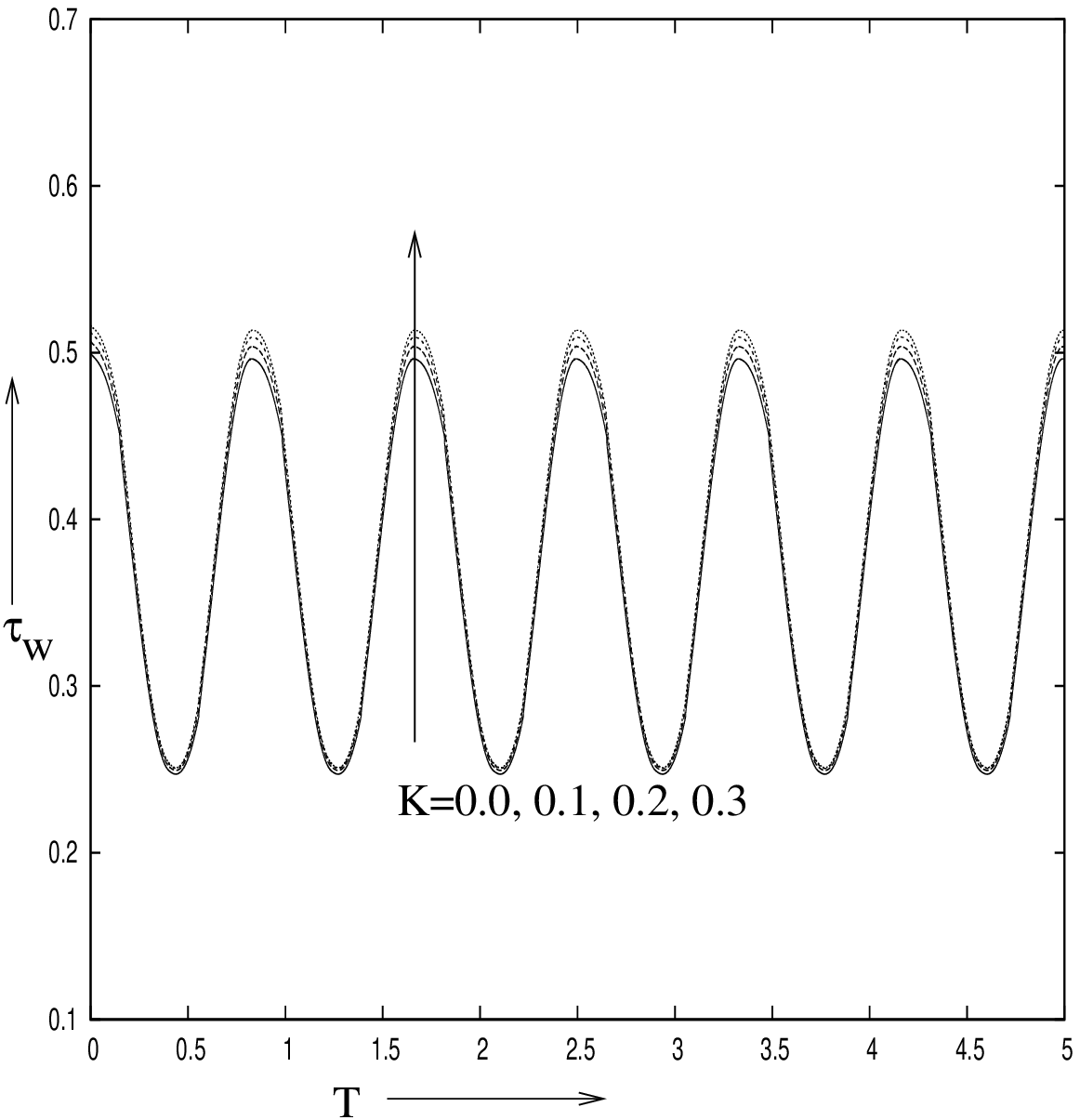}\\
Fig. 5 Variation of wall shear stress $\tau_{w}$ (a) along the
axis of the artery for different height of the stenosis $\delta$,
with  $H=2.0$, $a_{0}=1.0$, $K=0.1$; (b) with time $T$ for
different values of $H$, when $\delta=0.25$, $a_{0}=1.0$, $K=0.1$;
(c) with time $T$ for different amplitude of body acceleration
$a_{0}$, when $\delta=0.25$, $H=2.0$, $K=0.1$; (d) with time $T$
for different values of $K$, when $\delta=0.25$, $a_{0}=1.0$,
$H=2.0$.
\end{center}
\end{minipage}\vspace*{.5cm}\\
\newpage

\begin{minipage}{1.0\textwidth}
   \begin{center}
(a)\includegraphics[width=3.0in,height=2.0in ]{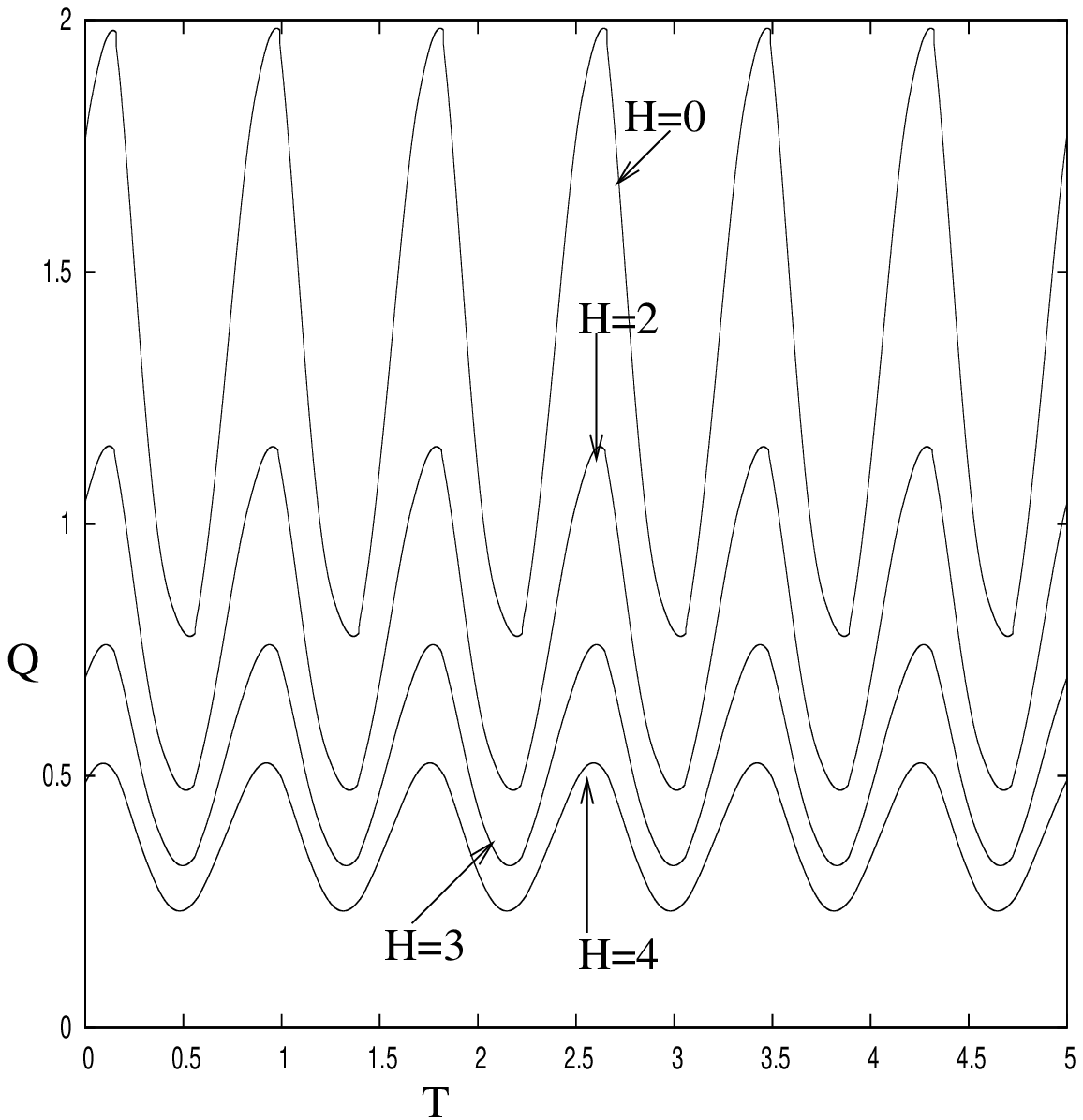}~~~
(b)\includegraphics[width=3.0in,height=2.0in ]{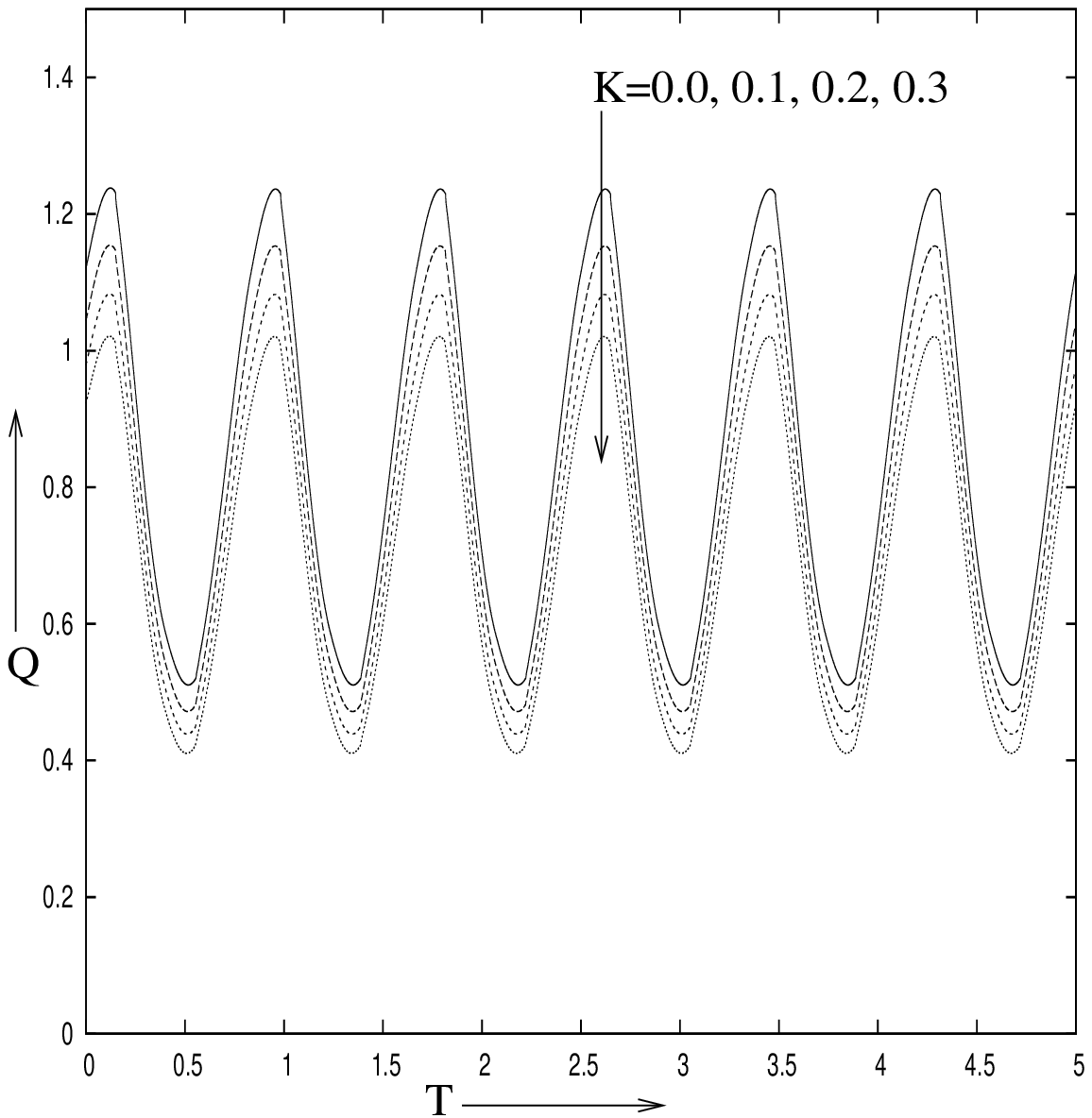}~~~\\
(c)\includegraphics[width=3.0in,height=2.0in ]{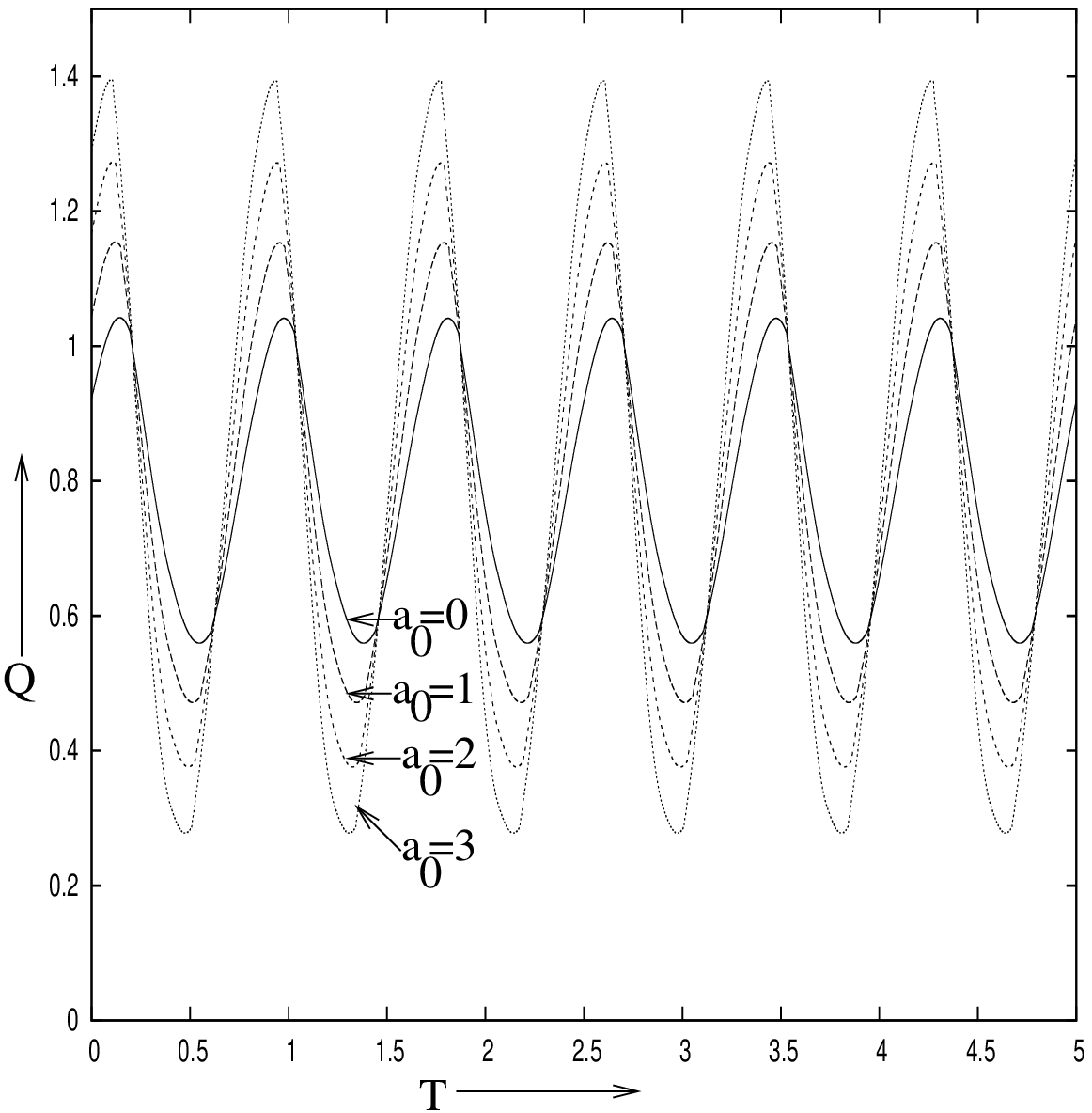}~~~~~\\
Figs. 6 Variation of volumetric flow rate $Q$ with time $T$(a) for
different values of $H$, when $K=0.1$, $a_{0}=1.0$; (b) for
different values of $H$, when $a_{0}=1.0$, $H=2.0$ and (c) for
different values of $a_{0}$, when $K=0.1$, $H=2.0$.
\end{center}
\end{minipage}\vspace*{.5cm}\\
\newpage

\begin{minipage}{1.0\textwidth}
(a)\includegraphics[width=3.0in,height=2.0in ]{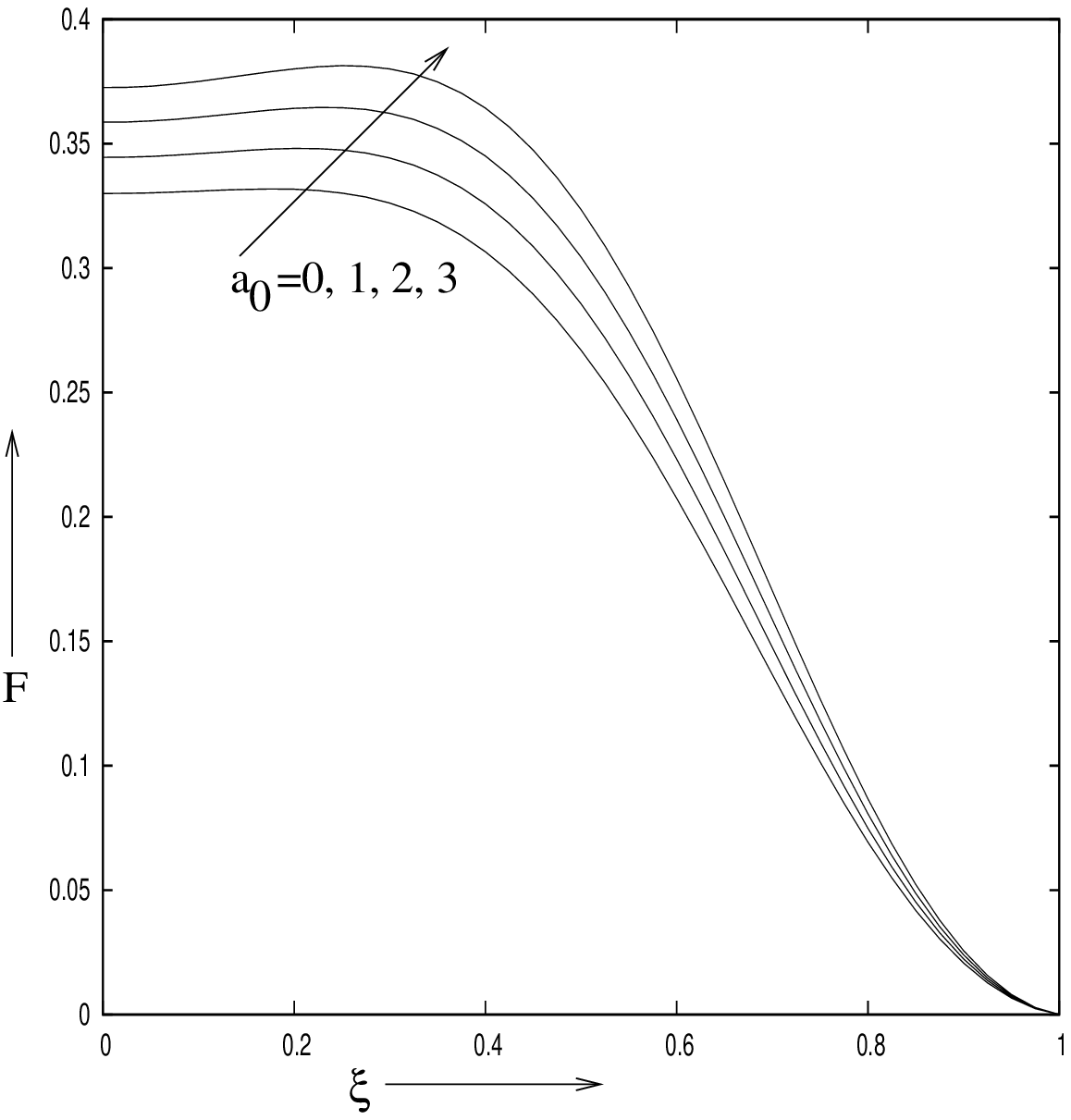}~~~
(b)\includegraphics[width=3.0in,height=2.0in ]{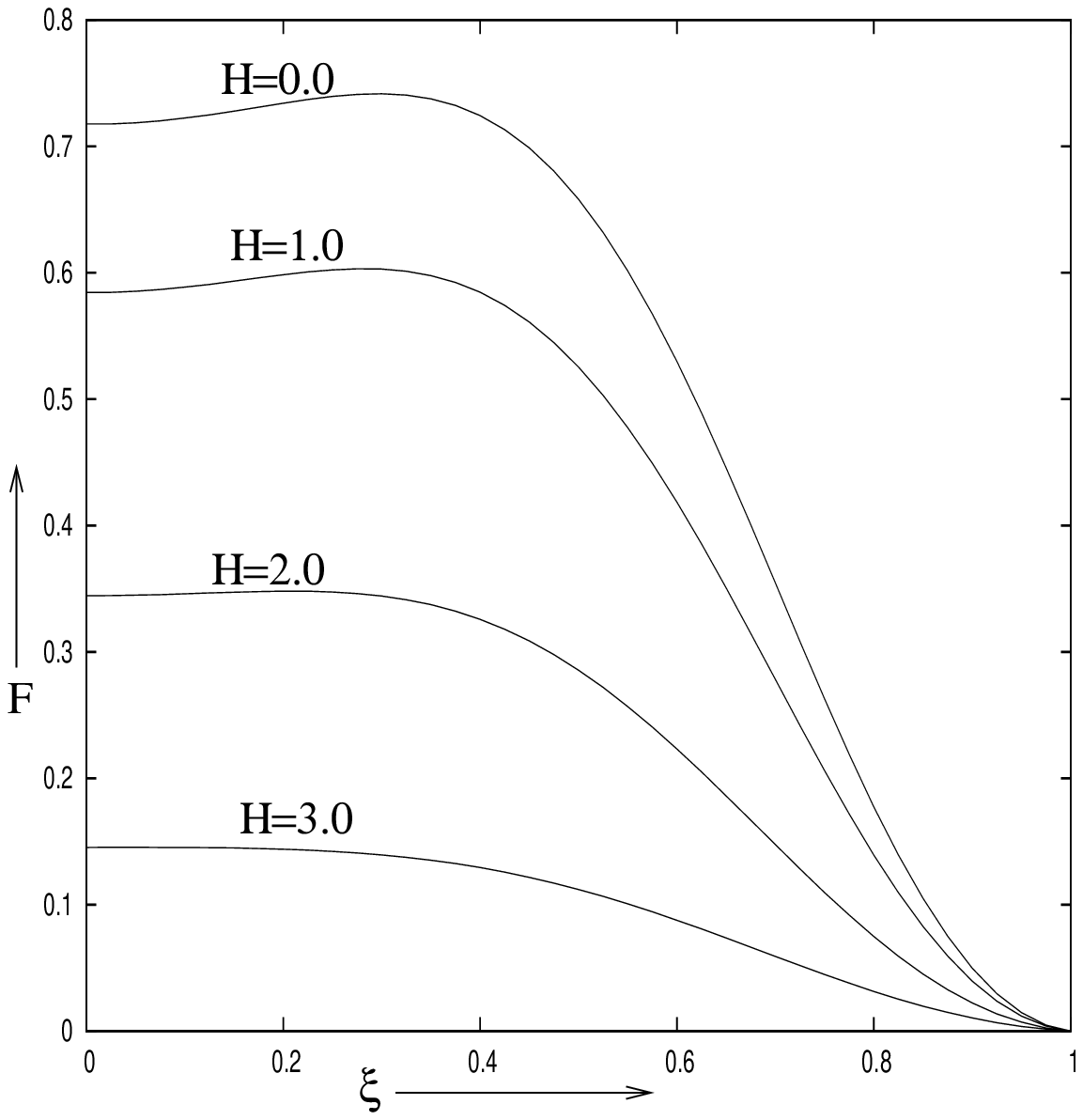}\\
Fig. 7 Variation of fluid acceleration $F$ along the radial
direction for (a) different amplitude of body acceleration $a_{0}$
when $H=2.0$ and (b) for different values of Hartmann number $H$
when $a_{0}=1.0$.
\end{minipage}\vspace*{.5cm}\\
\newpage

\begin{minipage}{1.0\textwidth}
   \begin{center}
(a)\includegraphics[width=3.0in,height=2.0in ]{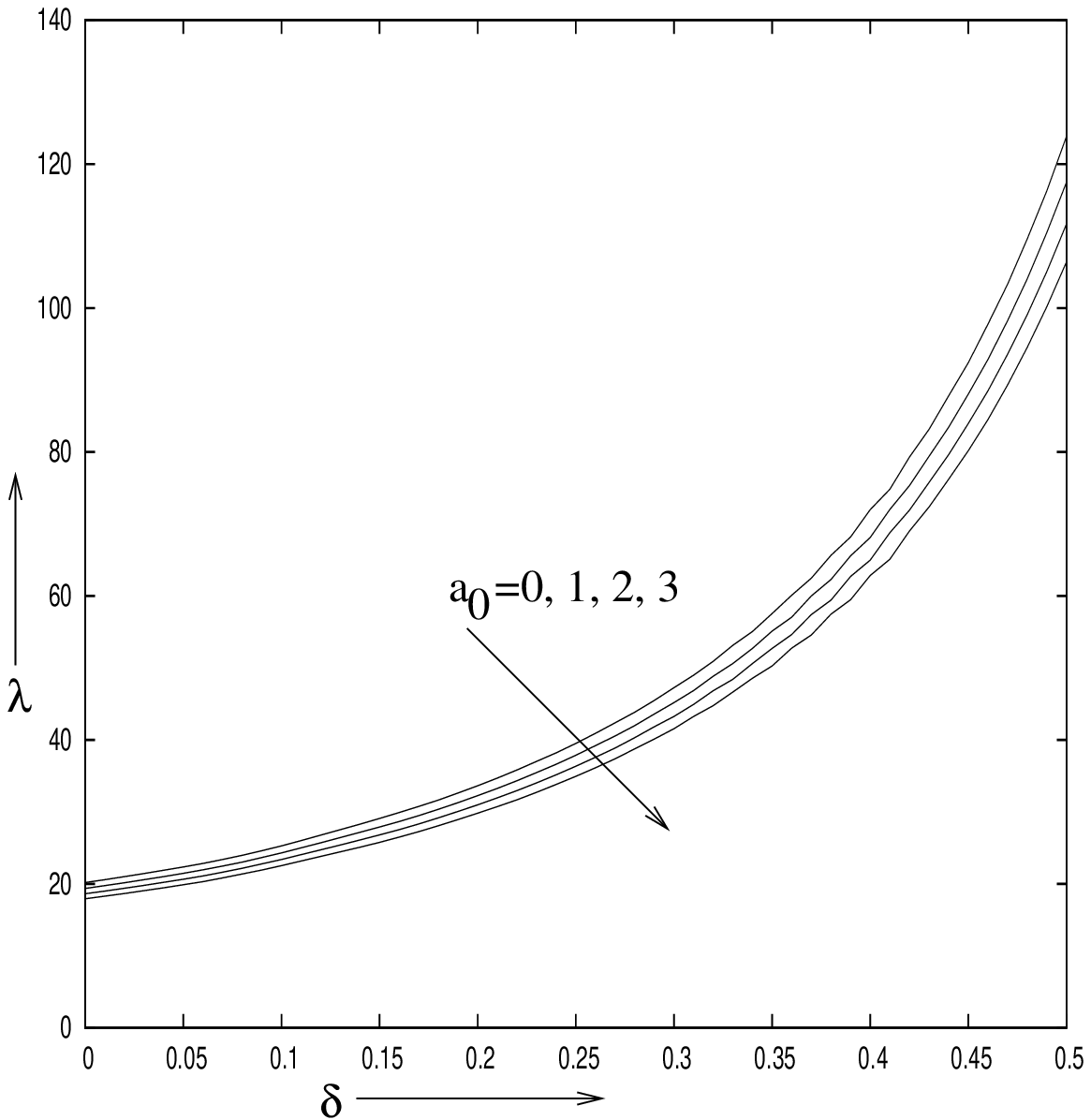}~~~
(b)\includegraphics[width=3.0in,height=2.0in ]{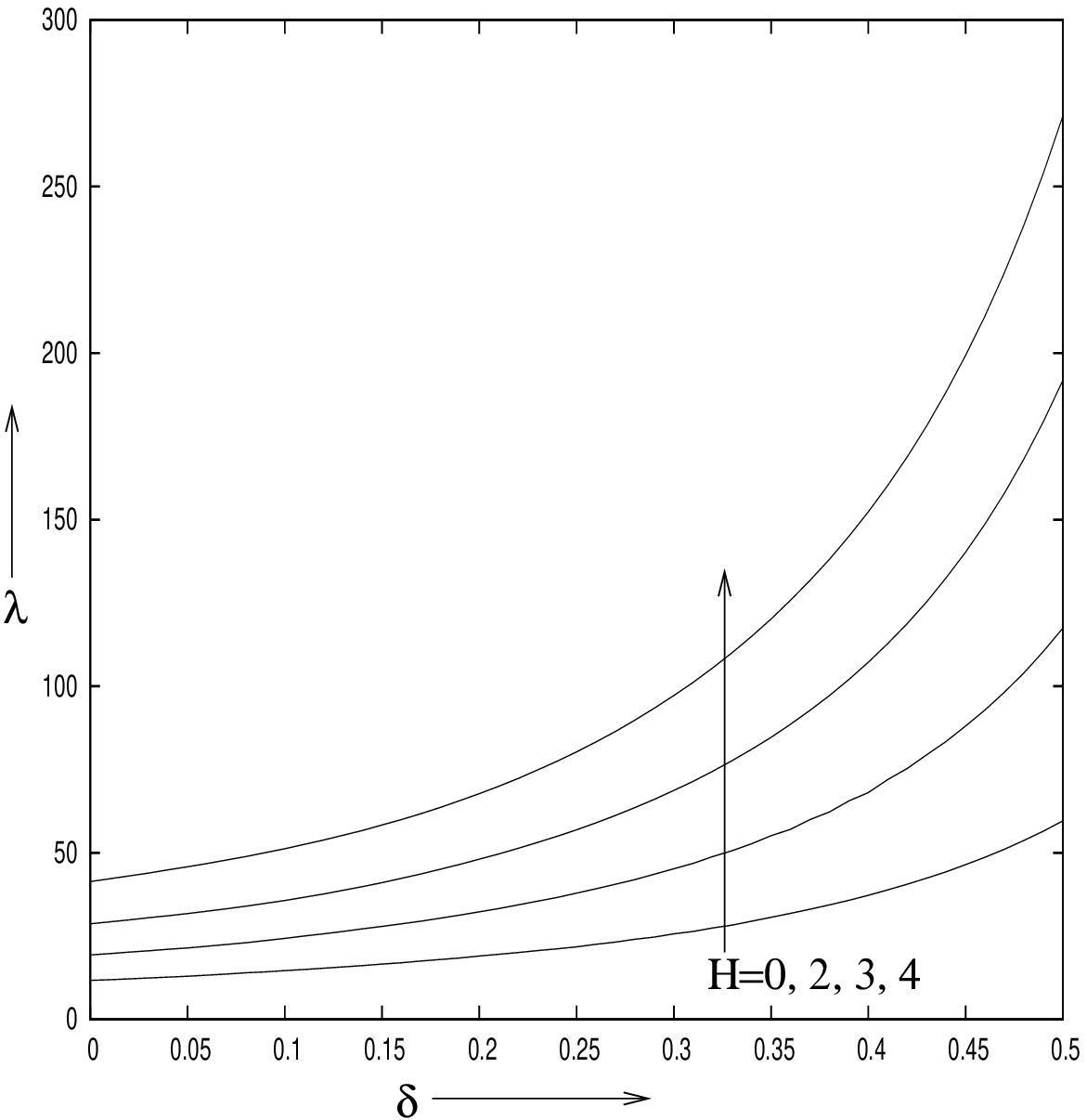}\\
(c)\includegraphics[width=3.0 in,height=2.0in ]{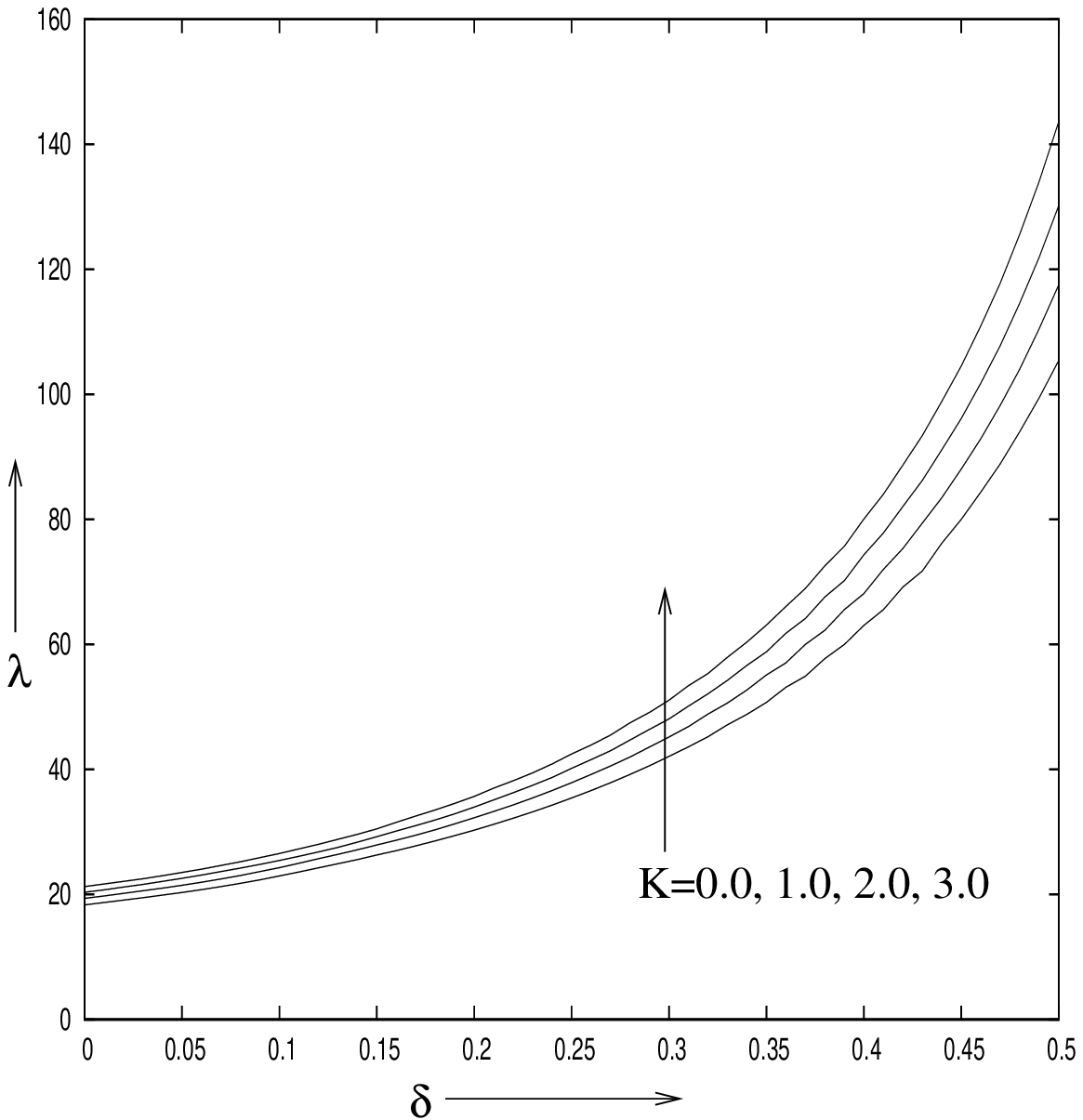}\\
Fig. 8 Variation of flow resistance $\lambda$ with depth of the
stenosis (a) for different values of $a_{0}$, when $H=2.0$,
$K=0.1$, (b) for different values of Hartmann number $H$, when
$a_{0}=1.0$, $K=0.1$ and (c) for different values of $K$, when
$H=2.0$, $a_{0}=1.0$.
\end{center}
\end{minipage}\vspace*{.5cm}\\

\end{document}